\documentclass[11pt,a4paper]{article}
\usepackage[english]{babel}
\usepackage{amsfonts,a4wide,bm,url}
\usepackage{amsmath}
\usepackage{amssymb}
\usepackage{amsthm}

\newtheorem{theorem}{Theorem}[section]

\newtheorem{corollary}[theorem]{Corollary}
\newtheorem{lemma}[theorem]{Lemma}
\newtheorem{proposition}[theorem]{Proposition}
\newtheorem{definition}[theorem]{Definition}
\newtheorem{remark}[theorem]{Remark}

\newcounter{c}
\newcounter{tmpabcd}

\newcounter{tmpnum}

\newcounter{tmprome}

\newcommand{\T}{\ensuremath{\mathcal{T}}}
\newcommand{\Tbz}{\T_{\b{z}}}

\newcommand{\car}{{\rm char\,}}

\newcommand{\mcc}{\mathbb{C}}

\newcommand{\mpp}{\mathbb{P}}

\newcommand{\mnn}{\mathbb{N}}

\newcommand{\mrr}{\mathbb{R}}

\newcommand{\mzz}{\mathbb{Z}}

\newcommand{\mqq}{\mathbb{Q}}

\newcommand{\diff}[1]{\frac{\partial}{\partial #1}}

\newcommand{\eqdef}{\ensuremath{\stackrel{\mathrm{def}}{=}}}

\newcommand{\rk}{\operatorname{rk}}

\newcommand{\Dist}{\ensuremath{{\rm Dist}}}
\newcommand{\dist}{\ensuremath{{\rm dist}}}

\newcommand{\dd}{\ensuremath{{\rm deg}}}

\renewcommand{\b}[1]{{{#1}}}

\newcommand{\hidden}[1]{}

\newcommand{\e}{{\rm e}}

\newcommand{\V}{\ensuremath{\mathcal{V}}}
\newcommand{\I}{\ensuremath{\mathcal{I}}}

\newcommand{\A}{\ensuremath{\mathcal{A}}}
\newcommand{\AnneauDePolynomes}{\A}

\newcommand{\Z}{\ensuremath{\mathcal{Z}}}

\newcommand{\Cr}{C_{\mbox{reg}}}
\newcommand{\Ciso}{C_{\mbox{iso}}}

\newcommand{\kk}{\ensuremath{\Bbbk}}

\newcommand{\KK}{\ensuremath{\overline{\kk((\b{z}))}}}

\newcommand{\idp}{\ensuremath{\mathcal{P}}}
\newcommand{\idq}{\ensuremath{\mathcal{Q}}}

\newcommand{\Spec}{\ensuremath{{\rm Spec}}}
\newcommand{\Ass}{\ensuremath{{\rm Ass}}}

\newcommand{\codim}{\ensuremath{{\rm codim}}}
\newcommand{\ul}[1]{\underline{#1}}
\newcommand{\ull}[1]{\underline{\b{#1}}}
\newcommand{\ullt}[1]{\underline{\b{#1}}}
\newcommand{\ol}[1]{\overline{#1}}

\newcommand{\bt}[1]{{\b{\tilde{#1}}}}

\newcommand{\ord}{\ensuremath{{\rm ord}}}
\newcommand{\ordz}{\ensuremath{{\rm ord_{\b{z}=0}}}}
\newcommand{\Ord}{\ensuremath{{\rm Ord}}}
\newcommand{\Ordz}{\ensuremath{{\rm Ord_{z=0}}}}
\newcommand{\trdeg}{\ensuremath{{\rm tr.deg.}}}

\newcommand{\Talg}{\ensuremath{\ul{\T}^*}}
\newcommand{\irrT}{\ensuremath{\irr \, \T}}

\newcommand{\irr}{{\rm irr}}
\newcommand{\rg}{{\rm rk}}

\newcommand{\eq}{{\rm eq}}

\begin{document}

\markboth{Evgeniy Zorin}
{Zero order estimates for analytic functions}

\title{ZERO ORDER ESTIMATES FOR ANALYTIC FUNCTIONS
}

\author{EVGENIY ZORIN\footnote{Institut de math\'ematiques de Jussieu, Universit\'e Paris 7, Paris, France.} \footnote{E-mail: evgeniyzorin@yandex.ru}} 
\date{}
\maketitle




\maketitle

\begin{abstract}
The primary goal of this paper is to provide a general multiplicity estimate. Our main theorem allows to reduce
a proof of multiplicity lemma to the study of ideals stable under some appropriate transformation of a polynomial ring.
In particular, this result leads to a new link between the theory of polarized algebraic dynamical systems and transcendental number theory.
On the other hand, it allows to establish an improvement of Nesterenko's conditional result on solutions of systems of differential equations. We also deduce, under some condition on stable varieties, the optimal multiplicity estimate in the case of generalized Mahler's functional equations, previously studied by Mahler, Nishioka, T\"opfer and others. Further, analyzing stable ideals we prove the unconditional optimal result in the case of linear functional systems of generalized Mahler's type. The latter result generalizes a famous theorem of Nishioka (1986) previously conjectured by Mahler (1969), and simultaneously it gives a counterpart in the case of functional systems for an important unconditional result of Nesterenko (1977) concerning linear differential systems.

In summary, we provide a new universal tool for transcendental number theory, applicable with fields of any characteristic. It opens the way to new results on algebraic independence, as shown in Zorin (2010).
\end{abstract}




\section{Introduction}

\subsection{General Background}

Classical methods of transcendental number theory estimate the transcendence degree of a field extension $\kk(f_1(\alpha),...,f_n(\alpha))/\kk$, where $\kk$ denotes a base field (e.g. a number field), $f_1,...,f_n\in \kk[[z]]$ are analytic functions
and $\alpha \in \overline{\kk}$ is algebraic over $\kk$ (or more generally, sometimes the same problem is studied for $\alpha\in L$, where $L$ is a field containing~$\kk$). These methods also provide measures of algebraic independence, i.e. lower bounds for
$|P(f_1(\alpha),...,f_n(\alpha))|$ (where $P \in \kk[X_1,\dots,X_n]$ denotes a non-zero polynomial)
in terms of degree and height of $P$. In general, these methods need a zero order estimate (or \emph{multiplicity estimates}, \emph{multiplicity lemmas}) for auxiliary function, which is usually of the form $P(f_1(z),...,f_n(z))$ (where $P$ denotes again a polynomial in $n$ variables). We say that one has a \emph{multiplicity lemma} for functions $f_1(z),\dots,f_n(z)$ if the following estimate holds for any non-zero polynomial $P\in \kk[X_1,\dots,X_n]$:
\begin{equation} \label{intro_ordzP}
    \ordz P(f_1(z),...,f_n(z))\leq f(\deg P),
\end{equation}
where $f:\mnn\rightarrow\mrr^+$ denotes a function depending on degree of $P$. One often introduces a more sophisticated control: for example function $f$ in the r.h.s. of~(\ref{intro_ordzP}) may depend on a \emph{height} of polynomial $P$ as well; or one may separate variables $X_1,\dots,X_n$ in two or more groups and introduce in the r.h.s. of~(\ref{intro_ordzP}) a function depending on the corresponding partial degrees.

There even exists a general method that uses this type of estimates  in order to establish transcendence or algebraic independence results (cf.~\cite{PP1997} and~\cite{PP_KF}).

Within this paper we prove a new multiplicity lemma applicable in a pretty general situation (Theorem~\ref{LMGP}; a slightly simplified version is given in Theorem~\ref{LMGP_simplified}). In particular, our multiplicity lemma answers the question raised in \cite{PP_KF} (see discussion at the end of~\S3 \emph{loc.cit.}) and complements in many interesting situations the method presented in this reference. Furthermore, this theorem allows to deal with all the classical cases mentioned below.

Our general theorem reduces the proof of multiplicity estimates to the study of ideals stable under some appropriate transformation of a polynomial ring. It generalizes a famous result of Yu.Nesterenko on differential systems satisfying a so called \emph{$D$-property}, see Chapter~10 in \cite{NP}. At the same time our theorem creates a new link between polarized algebraic dynamical systems and the theory of transcendental numbers (see Remark~\ref{rem_link_DA_PADS}).

To begin with, we give a brief account of known results postponing a more detailed presentation of our results to the second part of the introduction (see subsection~\ref{intro_section_presentation}).

Zero order estimates started to play an important role in works of C.L.Siegel~\cite{Si1929, Si1932}, A.O. Gelfond and Th.Schneider~\cite{Sch1934}. These authors were especially interested in zero order estimates for exponential polynomials. This work was completed by R.Tijdeman~\cite{Tijd1973}, whose well known theorem on the number of zeros of exponential polynomials generalizes the results of Gelfond.

 In~1954 A.B.Schidlovskii proved the first general estimate for the number of zeros of linear forms in solutions of a linear differential system. Then Yu.V.Nesterenko generalized his result from linear forms to polynomials introducing algebraic methods~\cite{N1972, N1973, N1974, N1977} in this context. In particular, he introduced the crucial notion of a stable ideal linking it to a differential ideal in Kolchin's theory.

In~1980 W.D.Brownawell and D.W.Masser~\cite{Brownawell1982,BrownawellMasserI1980,BrownawellMasserII1980}
gave the estimate
\begin{equation} \label{introduction_BrMas_est}
\ord_{z=0} P(f_1(z),...,f_n(z))\leq C\left(\deg P+1\right)^{2^n}
\end{equation}
for analytic (at $z=0$) functions $f_1,\dots,f_n$ satisfying a system of polynomial differential equations
\begin{equation} \label{BMsystem}
    f'_i=F_i(f_1,\dots,f_n),\quad i=1,\dots,n,
\end{equation}
where $F_i\in\mcc[X_1,\dots,X_n]$.
The estimate~(\ref{introduction_BrMas_est}) was improved to $\leq C\left(\deg P+1\right)^n$ by Nesterenko~\cite{N1988}. These results are based on the study of ideals in polynomial rings.

Yu.Nesterenko~\cite{N1996} also proved the estimate
\begin{equation} \label{bestSiegel}
\ord_{z=0} P(z,f_1(z),...,f_n(z)) \leq C\left(\deg_z P+1\right)\left( \deg_{\ul{f}} P + 1 \right)^n
\end{equation}
(optimal up to a multiplicative constant) for solutions $\ul{f}$
of differential system
\begin{equation} \label{relsNesternko}
{\frac{\partial f_i}{\partial z} =
\frac{A_i(z,f_1(z),...,f_n(z))}{z,A_0(f_1(z),...,f_n(z))}},
\end{equation}
$i=1,...,n$, where $A_0,...,A_n$ denote polynomials, under the condition that some
differential operator $D$ associated to~(\ref{relsNesternko})
has no stable ideals with a big order of annulation along $\ul{f}$ in
$z=0$ (see~\cite{NP}, chapter~10 or Theorem~\ref{theoNesterenko_classique} in subsection~\ref{subsection_ApplicationsDifferential} hereafter). This result allowed him to establish the following theorem: for all $q\in\mcc$, $0<|q|<1$ one has $\trdeg_{\mqq}\mqq(q,R(q),Q(q)R(q))\geq 3$, where $P$, $Q$ et $R$ denote Ramanujan functions (cf.~\cite{N1996}; one can also find the deduction of this result with a general method in~\cite{PP_KF}, this time admitting a multiplicity lemma). In particular, substituting $q=\e^{-2\pi}$ one obtains the algebraic independence over $\mqq$ of numbers $\e$, $\e^{\pi}$ and $\Gamma(\frac14)$ (see the beginning of Chapter~3, \cite{NP}).

If $A_0(0,f_1(0),...,f_n(0))\ne 0$ in~(\ref{relsNesternko}), one readily verifies the necessary condition on $D$-stable ideals deducing an unconditional multiplicity estimate (see Example~1 in Chapter~10 \cite{NP}). Also Nesterenko established this condition in the case when all the polynomials $A_i$ in~(\ref{relsNesternko}) are linear in $\ul{f}$~\cite{N1974}. For these purposes he used Galois theory of Picard-Vessiot extensions~\cite{Kol1948,Kol1973}, a result from the domain of differential algebra.

On the other hand, in the framework of Mahler's method K.Nishioka proved~\cite{Ni1986} the estimate~(\ref{bestSiegel}) for
$\ul{f}\in\kk[[{z}]]^n$ (where $\kk$ denotes a field of characteristic~0) satisfying a system of linear functional equations
\begin{equation} \label{relsNishioka}
{f_i}({z^d}) = \sum_{j=1}^n a_{ij}{{f_j}(z)}
\end{equation}
($i=1,...,n$), where $a_{ij}\in\kk(z)$, $d \geq 2$ an integer. This result was previously conjectured by K.Mahler~\cite{Mah1969}.

More generally, within the framework of Mahler's method one considers the following system of functional equations
\begin{equation} \label{relsTopfer}
{f_i}({p(z)}) =
\frac{A_i(z, f_1(z),...,f_n(z))}{A_0(z, f_1(z),...,f_n(z))},\index{Systeme d@Syst\`eme d' !equations de Mahler@\'equations de Mahler}
\end{equation}
$i=1,\dots,n$, where $A_0,...,A_n$ denote polynomials with $\deg_z A_i\leq s$, $\deg_{\ul{X}}A_i\leq t$ and $p(z)$ denotes a rational fraction, $\delta=\ord_{z=0} p({z})\geq 2$  and $d=\deg p(z)$. 

K.Nishioka proved~\cite{Ni1990}, for $p(z)=z^d$ (where $d\geq 2$ is supposed to be an integer) and $t^n<d$, that the following estimate holds for all non-zero polynomials $Q(z,X_1,...,X_n)$ satisfying $\deg_z Q\leq M$, $\deg_{\ul{X}}Q\leq N$ (where $M\geq N\geq 1$):
\begin{equation} \label{intro_estNishioka}
    \ord_{z=0} Q(z,f_1(z),...,f_n(z))\leq c_0MN^{n\log d/(\log d-n\log t)}.
\end{equation}

For a general rational substitution $p(z)$ Th.T\"opfer in~1998 proved~\cite{ThTopfer} the estimate
\begin{equation} \label{intro_estTopfer}
    \ord_{z=0} Q(z,f_1(z),...,f_n(z))\leq c_0MN^{n\log d/(\log\delta-n\log t)},
\end{equation}
where $d\eqdef\deg p$ and $\delta\eqdef\ord_{z=0}p$ (this time assuming $t^n<\delta$). Note that for $d=\delta$ (in particular, if $p(z)=z^d$) estimate~(\ref{intro_estTopfer}) coincides with~(\ref{intro_estNishioka}).

The results of Nishioka and T\"opfer are not the best possible, except for the case $d=\delta$ and $t=1$ (e.g. this is the case of relations~(\ref{relsNishioka})). In the other cases these estimates have the exponent strictly larger than $n$, so one can hope to ameliorate it.

All these results are obtained in characteristic~$0$. In a similar vein other authors study the results of the same types in the case of a strictly positive characteristic: P.-G.Becker~\cite{B1994}, V.Bosser and F.Pellarin~\cite{BP2008, BP2009}, F.Pellarin~\cite{Pellarin2009, Pellarin2010}.

Weakening the assumption $\ord_{z=0}p(z)\geq 2$ in~(\ref{relsTopfer}) to $\ord_{z=0}p(z)\geq 1$ we obtain another important case. Choosing $p(z)=qz$ (where $q\in\kk\setminus\{0\}$) and polynomials $A_i$ with $\deg_{\ul{X}}A_i=1$, we reduce~(\ref{relsTopfer}) to a system of equations in $q$-differences. This type of relations is also largely studied~\cite{AV2003,ATV2007,Bertrand2007}.

Recently P.Philippon proved a transference lemma~\cite{PP}, i.e. a statement that allows to construct an algebraic curve having a large contact number with analytic germ $$z\rightarrow\left(z,f_1(z),\dots,f_n(z)\right),$$ provided that there is a polynomial $Q$ with $\ordz Q(z,f_1(z),\dots,f_n(z))$ big enough. We show how a statement of this type allows one to establish multiplicity lemmas in many situations, encompassing all the classical cases mentioned above. However, we shall face a subtle problem concerning varieties (or ideals) stable under the action of some associated transformation.

\subsection{Presentation of results} \label{intro_section_presentation}

Our main result is a general multiplicity lemma (Theorem~\ref{LMGP}; a simplified version is given in Theorem~\ref{LMGP_simplified}). It deals with transformations $\phi:\AnneauDePolynomes\rightarrow\AnneauDePolynomes$ (where $\AnneauDePolynomes$ denotes a polynomial ring $\kk[X'_0,X'_1,X_0,\dots,X_n]$ and $\kk$ is a field) of a pretty general type. We introduce a notion of transformations \emph{correct} with respect to some ideal $I$. This notion is introduced in Definition~\ref{defin_phiestcorrecte}, at this point just note that this class includes all differential operators (for all the ideals $I$; see Corollary~\ref{exemple_OpDiff_Correcte}) as well as algebraic morphisms (the latter under some mild restriction, see Corollary~\ref{exemple_ApBirr_Correcte}). For example, condition~(\ref{intro_theoLMGP_condition_de_correctitude}) in Theorem~\ref{LMGP_simplified} below is automatically satisfied if $\phi$ is a differential operator.
We impose a natural hypothesis on a transformation $\phi$: it  should not too rapidly increase the degree of polynomials (cf.~(\ref{degphiQleqdegQ})) and it
should not too rapidly decrease an order of annulation at $\ul{f}=(1,z,1,f_1,\dots,f_n)\in\kk[[z]]^{n+3}$ of bi-homogeneous polynomials from $\AnneauDePolynomes$ (cf.~(\ref{condition_T2_facile})).

A key notion of our formal multiplicity lemma is a \emph{$\phi$-stable ideal}.
We say that an ideal $I\subset\A$ is \emph{$\phi$-stable} if $\phi(I)\subset I$ (cf. Definition~\ref{definIdealTstable}).
In the same way, let $\mpp$ denotes a multiprojective space (for example, $\mpp=\mpp^n$ or $\mpp=\mpp^1\times\mpp^n$) and $\T:\mpp\rightarrow\mpp$ be a rational application. We say that a variety $V\subset\mpp$ is $\T$-stable, if Zariski closure of $\T(V)$ coincides with $V$ (see Remark~\ref{rem_TV_simplified} and Definition~\ref{definVarieteTstable}). Relations between the notion of $\phi$-stable ideal and that of $\T$-stable variety are discussed in Remark~\ref{varietestable_et_idealstable}.

We shall systematically use the quantity $\ord_{\ul{f}}\idq$, where $\idq$ denotes a bihomogeneous ideal of $\A$ and $\ul{f}=(1,z,1,f_1(z),\dots,f_n(z))\in\kk[[z]]^{n+3}$; note that $\ul{f}$ can be considered as a system of projective coordinates of a point at $\mpp^1_{\kk[[z]]}\times\mpp^n_{\kk[[z]]}$. This quantity is introduced in Definition~\ref{defin_ord_xy} (see also \cite{EZ}, chapter~1, \S~3 for more details). Here we just mention that
$$
\ord_{\ul{f}}\idq\leq\min_{P\in\idq}\ordz P(\ul{f}).
$$
In particular, this means that \emph{$D$-property} (see Chapter~3 of~\cite{NP} for definition; $D$ is supposed to be a differential operator) implies for any $D$-stable ideal $\ord_{\ul{f}}\idq\leq K$ for some constant $K$ (in fact for the one provided by the \emph{$D$-property}). Hence, $D$-property implies~(\ref{intro2_RelMinN2}) in Theorem~\ref{LMGP_simplified} below.

Here is a simplified statement of our main result (note that this version still allows one to deduce, for example, our improvement of Nesterenko's result, Theorem~\ref{LMGPD}). The full version with additional parameters of control 
is given in Theorem~\ref{LMGP}. We explain how to obtain Theorem~\ref{LMGP_simplified} from Theorem~\ref{LMGP} in Remark~\ref{rem_n1}.

\begin{theorem}[Formal multiplicity lemma, simplified version]\label{LMGP_simplified}
Let $\kk$ be a field and $\ul{f}=(1,z,1,f_1,\dots,f_n)\in\kk[[z]]^{n+3}$. Let $\phi:\A\rightarrow\A$ be a (set-theoretical) transformation of a polynomial ring $\A=\kk[X_0',X_1'][X_0,\dots,X_n]$. Assume that $\phi$ satisfies~(\ref{degphiQleqdegQ}) and~(\ref{condition_T2_facile}).

We shall consider the polynomial ring $\A$ as bi-graduated with respect to $\left(\deg_{\ul{X}'},\deg_{\ul{X}}\right)$. Let $C_0 \in \mrr^+$ be a constant
such that for all bi-homogeneous prime ideal $\idq \subset \AnneauDePolynomes$ of rank $n$ one has
\begin{equation} \label{intro_theoLMGP_condition_de_correctitude}
\ord_{\ull{f}}\idq \geq C_0 \Rightarrow \mbox{ the transformation $\phi$ is correct with respect to }\idq.
\end{equation}
Suppose that there exists a
constant $K_0 \in \mrr^{+}$ (depending on $\phi$ and $\ull{f}$ only) with the following property: for all equidimensional bi-homogeneous $\phi$-stable ideal $I\subset\AnneauDePolynomes$ such that all its associated prime ideals satisfy
\begin{equation} 
\ord_{\ull{f}}\idq \geq C_0,
\end{equation}
there exists a prime factor $\idq\in\Ass(\AnneauDePolynomes/I)$ that satisfies
\begin{equation} \label{intro2_RelMinN2}
\ord_{\ull{f}}(\idq) < K_0\left(\dd_{(0, n-\rg\idq+1)}(\idq)+\dd_{(1, n-\rg\idq)}(\idq)\right).
\end{equation}

Then there exists a constant $K>0$ such that for all $P \in
\AnneauDePolynomes \setminus \{0\}$ one has
\begin{equation} \label{intro_LdMpolynome2}
\ordz(P(\ull{f})) \leq K\left(\deg_{\ul{X}'}P+1)+\epsilon\deg_{\ul{X}}P\right)\times(\deg_{\ul{X}} P + 1)^n,
\end{equation}
where $\epsilon$ may denote either $0$ or $1$ in function of the properties satisfied by the transformation $\phi$. In particular,~(\ref{intro_LdMpolynome2}) with $\epsilon=1$ is true in any case, and in the case when $\phi$ is a differential operator we have a stronger estimate: (\ref{intro_LdMpolynome2}) with $\epsilon=0$.
\end{theorem}

In Remark~\ref{rem_n1} we explain how to obtain Theorem~\ref{LMGP_simplified} from Theorem~\ref{LMGP}.

We also provide specialized versions of Theorem~\ref{LMGP} (one can substitute Theorem~\ref{LMGP_simplified} for the sake of simplicity) in the cases when $\phi$ is a differential operator (Theorem~\ref{LMGPD}) or algebraic morphism (Theorem~\ref{LMGPF}). In both cases we achieve to restrict the assumption~(\ref{intro2_RelMinN2}) (or equivalently~(\ref{RelMinN2}) for the same property in Theorem~\ref{LMGP}) to prime ideals only (or irreducible varieties in Theorem~\ref{LMGPF}), whereas in general theorem (Theorem~\ref{LMGP} or Theorem~\ref{LMGP_simplified}) one \emph{a priori} has to assure~(\ref{intro2_RelMinN2}) for all $\phi$-stable ideals.

\begin{remark} \label{rem_link_DA_PADS}
In case when $\phi=D$ is a differential operator, the notion of $D$-stable ideal coincides with the one of \emph{differential ideal} in the sense of differential algebra (for example, see~\cite{Kap1, Kol1942, Kol1948} for introduction). Nowadays this is an extensive and highly developed area of algebra~\cite{Cout1995,Mag1994,Mar1996} (for sure, even giving a panorama of this area demands at least a separate article). This theory plays an essential role in the proof by Nesterenko of his multiplicity lemma for solutions of linear differential equations~\cite{N1974} (more precisely, he used Galois theory for Picard-Vessiot extensions~\cite{Kol1948,Kol1973}).

In the case when $\phi$ is an algebraic morphism we arrive (with Theorem~\ref{LMGPF}) to the study of varieties stable under certain morphism. As this varieties are emerged in a (bi-)projective space they are \emph{polarized} in a sense explicated in~\cite{Zhang2006}. Thus we have in this case a \emph{polarized algebraic dynamical system}. The opposite is also true as shown in~\cite{Fakhruddin2003}: all such dynamical systems can be emerged in a projective space of a sufficiently big dimension.
Study and classification of polarized algebraic dynamical systems is a vast and extensively developed area 
(see~\cite{Zhang2006} for an overview and further references of the subject).
\end{remark}

Finally, in section~\ref{sectionNishioka} we study $\T$-stable varieties, for $\T:\mpp^1\times\mpp^n \rightarrow\mpp^1\times\mpp^n$ defined by
\begin{multline} \label{intro_defT2}
(X_0':X_1',X_0:...:X_n)\rightarrow\Big(A_0'(X_0',X_1'):A_1'(X_0',X_1'),\\
\sum_{j=0}^n a_{0j}(\ul{X}')X_j:...:\sum_{j=0}^n a_{nj}(\ul{X}')X_j\Big),
\end{multline}
where $A_i',a_{ij}(\ul{X}')\in\kk[X_0',X_1']$, $i=0,1$ and $a_{kj}(\ul{X}')\in\kk[X_0',X_1']$, $k,j=0,\dots,n$ denote homogeneous polynomials of degrees in $\ul{X}'$ respectively $r$ and $s$. We assume that $A_0'(X_0',X_1')$ and $\sum_{j=0}^n a_{0j}(\ul{X}')X_j$ both are non-zero polynomials.

The transformation $\T$ defined in~(\ref{intro_defT2}) is in a sense associated to the following system of functional equations~(\ref{relsTopfer}) considered within a framework of Mahler's method:
\begin{equation} \label{intro_relsTopfer2}
    \left(\sum_{j=0}^n a_{0j}(\ul{X}'){f}_j(\b{z})\right)f_i(p(\b{z}))=\sum_{j=0}^n a_{ij}(\ul{X}'){f}_j(\b{z}), \quad i=1,...,n,
\end{equation}
where $p(\b{z})=\frac{A_1'(1,\b{z})}{A_0'(1,\b{z})}$ (see Remark~\ref{rem_Mutual_Association} and Definition~\ref{def_Mutual_Association} for a general definition).

Our study of $\T$-stable varieties allows to verify condition~(\ref{RelMinN}) of Theorem~\ref{LMGPF} (a counterpart of condition~(\ref{intro2_RelMinN2}) of Theorem~\ref{LMGP_simplified}) in this case. So we establish the following unconditional multiplicity lemma:
\begin{theorem} \label{theoNishioka}
Let $\kk$ be a field of an arbitrary characteristic and $\T:\mpp^1_{\kk}\times\mpp^n_{\kk}\rightarrow\mpp^1_{\kk}\times\mpp^n_{\kk}$ a transformation defined by~(\ref{intro_defT2}). Suppose that there is a solution $(1,f_1(\b{z}),...,f_n(\b{z}))$ in algebraically independent power series of the system of functional equations~(\ref{intro_relsTopfer2}), which is associated to $\T$. Assume also that
\begin{equation}\label{theoNishioka_lambda_pgq2}
    \lambda:=\ordz p(\b{z})\geq 2.
\end{equation}

Then there exists a constant $K_1$ such that for any non-zero polynomial $P\in\AnneauDePolynomes$ one has
\begin{equation} \label{intro_theoNishioka_conclusion}
    \ordz(P(\ull{f})) \leq K_1(\deg_{\ul{X'}}P + \deg_{\ul{X}}P + 1)(\deg_{\ul{X}} P + 1)^n.
\end{equation}
\end{theorem}
Theorem~\ref{theoNishioka} was previously proved by Nishioka in the particular case $p(\b{z})=\b{z}^d$, ($d\geq 2$)~\cite{Ni1990}. T\"opfer considered more general systems but his estimates for $\ord_{\ull{f}}P$ have the exponent strictly bigger than $n$ (and so they are not optimal)~\cite{ThTopfer}. All these theorems were established in characteristic~0 only, contrarily to our Theorem~\ref{theoNishioka}.

Multiplicity estimates proved here allow to obtain new results on algebraic independence and improve the measures of algebraic independence. In particular, they are useful in the framework of Mahler's method~\cite{Pellarin2010} where they allow~\cite{EZ} to construct new large family of $n$-tuples that are \emph{normal} in the sense introduced by G.Chudnovsky~\cite{Ch1980}. We plan to return to this issue in a subsequent publication.

\section{Framework, definitions and first properties}

\subsection{General framework} \label{subsection_general_framework}

Within this subsection we fix formally the general framework imposed in this article.

Let $\kk$ be a (commutative) algebraically closed field and $\A$ be a ring of polynomials with coefficients in $\kk$: $\A=\kk[X_0',X_1'][X_0,...,X_n]$. We have marked here two groups of variables, $\ul{X}'$ and $\ul{X}$, because it is convenient for us to consider the ring $\A$ as bigraduated with respect to $\deg_{\ul{X}'}$ and $\deg_{\ul{X}}$. Note that $\kk$ may have a positive characteristic as well as characteristic~0.

\begin{remark}
The assumption that field $\kk$ is algebraically closed is not too restrictive. Indeed, having some field $k$, not necessary algebraically closed, we can apply our theorems with the algebraic closure $\kk:=\ol{k}$ to obtain the conclusion of multiplicity lemma, say~(\ref{intro_theoNishioka_conclusion}). This conclusion implies immediately the same statement for $k$, because $k\subset\kk$.

So our unconditional result, Theorem~\ref{theoNishioka}, is valid for any base field, whether it is algebraically closed or not. As for our conditional results, such as Theorem~\ref{LMGP_simplified}, Theorem~\ref{LMGP}, Theorem~\ref{LMGPD} or Theorem~\ref{LMGPF}, the only thing that one has to add to apply them to arbitrary field $k$ is that the imposed conditions on stable ideals (e.g.~(\ref{intro2_RelMinN2}) for Theorem~\ref{LMGP_simplified}) have to be verified for all the stable ideals in the ring $\ol{k}[X_0',X_1'][X_0,...,X_n]$, and not just in $k[X_0',X_1'][X_0,...,X_n]$.
\end{remark}


We fix a transformation $\phi:\A\rightarrow\A$ such that for all
$Q\in\AnneauDePolynomes$ one has
\begin{equation} \label{degphiQleqdegQ}
\begin{aligned}
&\deg_{\ul{X}} \phi(Q) \leq \mu \deg_{\ul{X}} Q,\\
&\deg_{\ul{X}'} \phi(Q) \leq \nu_0 \deg_{\ul{X}'} Q + \nu_1 \deg_{\ul{X}} Q
\end{aligned}
\end{equation}
with some constants $\mu$, $\nu_0$ et $\nu_1$.

\begin{remark} \label{rem_phi_settheoretical}
We would like to emphasize that the word \emph{transformation} means in this article simply a set-theoretical application, so $\phi$ is defined everywhere on the set $\A$ and gets values in the same set $\A$, but we do not suppose \emph{a priori} that $\phi$ respects any additional structure defined on $\A$, for example that one of the polynomial ring. All the properties imposed \emph{a priori} on $\phi$ are~(\ref{degphiQleqdegQ}) and~(\ref{condition_T2_facile}) below.
\end{remark}

Our principal result, Theorem~\ref{LMGP}, is proved for \emph{transformations} satisfying these (quite mild) restrictions, as well as one additional condition described in Definition~\ref{defin_phiestcorrecte}.

\begin{remark}
We denote by $\phi^N$ the $N$-th iteration of the transformation $\phi$.
\end{remark}

We fix a point $\ul{f}=(1:\b{z},1:f_1(\b{z}):...:f_n(\b{z})) \in \mpp^1_{\kk[[\b{z}]]}\times\mpp^n_{\kk[[\b{z}]]}$
and assume that there exist two constants $\lambda>0$ and $K_{\lambda}\geq 0$ such that
\begin{equation} \label{condition_T2_facile}
    \ordz \phi(Q)(\ul{f}) \geq \lambda \, \ordz(Q(\ul{f})).
\end{equation}
for all polynomials $Q\in\AnneauDePolynomes$ satisfying $\ordz(Q(\ul{f}))\geq K_{\lambda}$.

Two typical examples of such transformation $\phi$ are differential operators and algebraic morphisms.

Using recurrence on the condition~(\ref{degphiQleqdegQ}) we readily prove the following lemma:

\begin{lemma} \label{majorationphinQ} Let $N$ be a positive integer and $Q \in \AnneauDePolynomes$. Then
\begin{eqnarray}
  \deg_{\ul{X}}\phi^N(Q) &\leq& \mu^N\deg_{\ul{X}}Q, \label{majorationdegXphinQ} \\
  \deg_{\ul{X'}}\phi^N(Q) &\leq& \nu_0^N\deg_{\ul{X'}}Q+\nu_1\left(\sum_{i=0}^{N-1}\nu_0^{N-i-1}\mu^i\right)\deg_{\ul{X}}Q. \label{majorationdegXprimephinQ}
\end{eqnarray}
\end{lemma}
The proof of Lemma~\ref{majorationphinQ} is easy and left to the reader, but see Lemma~2.1 in~\cite{EZ} for a hint.

\begin{remark}
We can consider $\phi$ as acting on $\kk[\b{z}][X_0:...:X_n]$ by setting
\begin{equation}
\b{\phi}(\b{Q})=\phi\left(X_0'^{\deg_{\b{z}}\b{Q}}\b{Q}(X_1'/X_0')(X_0:...:X_n)\right)\Bigg|_{(X_0':X_1')=(1,\b{z})}
\end{equation}
for all $Q\in\kk[\b{z}][X_0:...:X_n]$.

This application $\b{\phi}$ satisfies
\begin{equation} \label{hphiQleqhQ}
\begin{aligned}
\deg_{\ul{X}} \b{\phi}(\b{Q}) &\leq \mu \deg_{\ul{X}} \b{Q},\\
h(\b{\phi}(\b{Q}))=\deg_{\b{z}} \b{\phi}(\b{Q}) &\leq \nu_0 \deg_{\b{z}} \b{Q} + \nu_1\deg_{\ul{X}} \b{Q} \\&\leq \nu_0 h(\b{Q}) + \nu_1\deg\b{Q}.
\end{aligned}
\end{equation}
\end{remark}

\subsection{Definitions and properties related to commutative algebra} \label{definitions_comm_algebra}

\begin{definition}
Let $I\subset\AnneauDePolynomes$ be a bi-homogeneous ideal.  We denote by $\V(I)$ the sub-scheme of $\mpp^1\times\mpp^n$ defined by $I$.
Conversely, for any sub-scheme $V$ of $\mpp^1\times\mpp^n$ we denote $\I(V)$ the bi-homogeneous saturated ideal in $\AnneauDePolynomes$ that defines $V$.
\end{definition}

\begin{definition} \label{def_I}
Let $V$ be a $\kk$-linear subspase of $\A$ and $\idp\subset\A$ a prime ideal. We define
\begin{equation*}
    I(V,\idp) \eqdef (V\AnneauDePolynomes_{\idp})\cap\AnneauDePolynomes,
\end{equation*}
where $\AnneauDePolynomes_{\idp}$ denotes the localization of $\AnneauDePolynomes$ by $\idp$ and $V\AnneauDePolynomes_{\idp}$ denotes
the ideal generated by $\AnneauDePolynomes_{\idp}$ by the elements of $V$.
\end{definition}


\begin{remark}
    Let $I_V$ denotes the ideal generated in $\A$ by $V$. Then ideal $I(V,\idp)$ is the intersection of the primary components of $I_V$ contained in $\idp$.
\end{remark}

\begin{definition}\label{definIdealTstable}
We say that an ideal $I \subset \AnneauDePolynomes$ is \emph{$\phi$-stable} if $\phi(I) \subset I$.
\end{definition}

\begin{definition} \label{def_eqI} Let $I$ be an ideal and
\begin{equation}
I = \mathcal{Q}_1 \cap \dots \cap \mathcal{Q}_r \cap
\mathcal{Q}_{r+1} \cap ... \cap \mathcal{Q}_{s}
\end{equation}
be its primary decomposition, where
$\mathcal{Q}_1$,...,$\mathcal{Q}_r$ are the primary ideals associated to the ideals of minimal rank (i.e. of rank
$\rg(I)$) and $\mathcal{Q}_{r+1}$,...,$\mathcal{Q}_{s}$
correspond to the components of rank
strictly bigger than $\rg(I)$.

We denote by
\begin{equation}
 \eq(I) \eqdef \mathcal{Q}_1 \cap \dots \cap \mathcal{Q}_r
\end{equation}
the equidimensional part of minimal rank of $I$.
\end{definition}

As it was mentioned above, principally our results deals with transformations satisfying one more special condition:

\begin{definition} \label{defin_phiestcorrecte} We say that an application $\phi:\A\rightarrow\A$ is {\it correct with respect to the ideal $\idp\subset\AnneauDePolynomes$} if for every ideal
$I$, such that all its associated primes are contained in $\idp$, the inclusion
\begin{equation}\label{defin_phi_phiI_subset_eqI}
    \phi(I)\subset\eq(I)
\end{equation}
implies
\begin{equation}\label{defin_phi_phieqI_subset_eqI}
    \phi(\eq(I))\subset\eq(I)
\end{equation}
(recall that $\eq(I)$ is introduced in Definition~\ref{def_eqI}).
\end{definition}

The restriction to correct transformations (with respect to some ideal $\idp$) is not too tough: as examples of correct transformations we can mention at least differential operators and (dominant) algebraic morphisms. In order to prove this we need the following technical lemma.

\begin{lemma} \label{lem_phicorrecte} Let $\idp\subset \AnneauDePolynomes$ be an ideal and $\phi:\AnneauDePolynomes \rightarrow \AnneauDePolynomes$
any transformation. Assume that for any equidimensional ideal $J$ such that all its associated primes are contained in $\idp$
\emph{one} of two following hypothesis is satisfied for all $x\in\idp$:
\begin{eqnarray}
x\phi(J) &\subset& \phi(x J)+J, \label{phicorrectetype1} \\ \text{ \emph{or} } \nonumber \\
\phi(x)\phi(J) &\subset& \phi(x J)+J \label{phicorrectetype2} \\ \mbox{ and } \rg \, \phi^{-1}(\idq)&=&\rg \idq \text{ for all } \idq\in\Ass(\AnneauDePolynomes/J). \label{phicorrectetype2_2}
\end{eqnarray}
Then, the transformation $\phi$ is correct with respect to $\idp$.
\end{lemma}
\begin{proof}
Let $I\subset\AnneauDePolynomes$ be an ideal such that all its associated primes are contained in $\idp$
and such that one has~(\ref{defin_phi_phiI_subset_eqI}).
In order to prove the lemma we have to verify~(\ref{defin_phi_phieqI_subset_eqI}).

Assume at first that $\phi$ satisfies~(\ref{phicorrectetype1}).
In this case we choose
\begin{equation} \label{x_equalize1}
x \in I:_{\A}\eq(I)
\end{equation}
such that
\begin{equation} \label{x_ni_AssJe}
x \mbox{ is a non-zerodivisor in } \AnneauDePolynomes/\eq(I)
\end{equation}
(a choice of one such $x$ is possible because the primes associated to
$I:_{\A}\eq(I)$ are of rank strictly bigger than that of primes associated to $\eq(I)$).

In view of the hypothesis~(\ref{phicorrectetype1}) (applied with $J=\eq(I)$)
\begin{equation*}
x \phi(\eq(I)) \subset \phi(x\,\eq(I))+\eq(I)
\end{equation*}
and with~(\ref{x_equalize1}), then~(\ref{defin_phi_phiI_subset_eqI}), we deduce
\begin{equation} \label{xphiJe_subset_Je1}
x\phi(\eq(I))\subset\phi(I)+\eq(I)\subset\eq(I)+\eq(I)\subset\eq(I).
\end{equation}
As we have chosen $x$ satisfying~(\ref{x_ni_AssJe}), we deduce from~(\ref{xphiJe_subset_Je1})
\begin{equation*} 
\phi(\eq(I))\subset\eq(I),
\end{equation*}
establishing~(\ref{defin_phi_phieqI_subset_eqI}) in this case (i.e. under the hypothesis~(\ref{phicorrectetype1})).

Let's consider now the second possibility from the statement: assume that we have hypothesis~(\ref{phicorrectetype2}). In this case we choose an
$x \in I:_{\A}\eq(I)$ such that
\begin{equation} \label{Tx_nest_divzero}
\phi(x) \mbox{ is a non-zerodivisor in } \AnneauDePolynomes/\eq(I)
\end{equation}
(the latter condition is equivalent to the fact $x \not \in {\phi}^{-1}(\idq)$ for all $\idq \in \Ass(\eq(I))$).
A choice of one such $x$ is possible because the primes associated to
$I:_{\A}\eq(I)$ are of rank strictly bigger than
that of primes associated to $\eq(I)$ and the hypothesis~(\ref{phicorrectetype2_2}) assures us that $\rg\left(\phi^{-1}(\idq)\right)=\rg \idq$
for all ideal $\idq \subset \idp$.

With our choice of $x$ we have $x\,\eq(I)\subset I$, hence we deduce the hypothesis~(\ref{phicorrectetype2}) (applied with $J=\eq(I)$)
\begin{equation*}
 \phi(x)\phi(\eq(I))\subset\phi(x\,\eq(I))+\eq(I)\subset\phi(I)+\eq(I)\subset\eq(I)+\eq(I)\subset\eq(I).
\end{equation*}
We conclude $\phi(\eq(I))\subset\eq(I)$, because in view of~(\ref{Tx_nest_divzero}) we have that $\phi(x)$ is not a zero divisor in $\AnneauDePolynomes/\eq(I)$.

Thus we have established~(\ref{defin_phi_phieqI_subset_eqI}) using the hypothesis~(\ref{phicorrectetype2}) and this completes the proof of Lemma~\ref{lem_phicorrecte}.
\end{proof}

In two following corollaries we deduce with Lemma~\ref{lem_phicorrecte} the statements that differential operators and algebraic transformations associated to rational applications (in the latter case under a mild restriction, see Corollary~\ref{exemple_ApBirr_Correcte}) are correct transformations.

\begin{corollary} \label{exemple_OpDiff_Correcte} Let $D:\AnneauDePolynomes \rightarrow \AnneauDePolynomes$ be a differential operator.
Then for any ideal $\idp\subset\AnneauDePolynomes$, the application $\phi=D$ is correct with respect to $\idp$.
\end{corollary}
\begin{proof}
Let $J\subset\AnneauDePolynomes$ be an ideal and $x\in\AnneauDePolynomes$. Using Newton-Leibnitz equality, one has for all $a\in J$
\begin{equation*}
    \phi(x\cdot a)=x\cdot\phi(a)+\phi(x)\cdot a.
\end{equation*}
As $J$ is an ideal, we obtain $\phi(x)\cdot a\in J$, and so
\begin{equation*}
    x\cdot\phi(a)=\phi(x\cdot a)-\phi(x)\cdot a\in \phi(xJ)+J.
\end{equation*}
Thus we have verified the condition~(\ref{phicorrectetype1}) of Lemma~\ref{lem_phicorrecte} and with this lemma we conclude that $\phi=D$
is correct with respect to any ideal $\idp\subset\AnneauDePolynomes$.
\end{proof}

\begin{corollary} \label{exemple_ApBirr_Correcte} Let $\T: \mpp^1\times\mpp^n \rightarrow \mpp^1\times\mpp^n$ be a rational application,
such that for a point $\ul{\alpha} \in \mpp^1\times\mpp^n$ one has the following property: each irreducible variety $V$ passing by $\ul{\alpha}$
satisfies
\begin{equation}\label{exemple_ApBirr_Correcte_dimVeqdimTV}
    \dim\T(V)=\dim V.
\end{equation}
We denote by $\idp$ the ideal of $\AnneauDePolynomes$
defining the point $\ul{\alpha}$. Then, the transformation $\phi=\T^*$ is correct with respect to $\idp$.
\end{corollary}
\begin{proof}
As $\Talg$ is an endomorphism of ring $\AnneauDePolynomes$, we immediately obtain (for any ideal $J\subset\AnneauDePolynomes$)
\begin{equation*}
    \phi(x)\phi(J)=\phi(xJ),
\end{equation*}
so the hypothesis~(\ref{phicorrectetype2})
of Lemma~\ref{lem_phicorrecte} is satisfied. For prime ideals we have $\V({\Talg}^{-1}(\idq))=\T(\V(\idq))$ and so
the hypothesis~(\ref{exemple_ApBirr_Correcte_dimVeqdimTV})
implies the hypothesis~(\ref{phicorrectetype2_2}) of Lemma~\ref{lem_phicorrecte}.
By this lemma we obtain that the transformation $\phi=\Talg$ is correct
with respect to the ideal $\idp$ from the statement.
\end{proof}

\vbox{
\begin{definition}\label{definDePP}\begin{enumerate}
  \item Let $\idp$ be a prime ideal of the ring $\AnneauDePolynomes$, $V$ a $\kk$-vectorial subspace
of $\AnneauDePolynomes$ and $\phi$ a transformation of the ring $\AnneauDePolynomes$ to itself. Then
\begin{equation} \label{definDePP_defin_e}
 e_{\phi}(V,\idp) \eqdef \max(e\,\vline\,\rg\left((V+\phi(V)+...+{\phi}^e(V))\AnneauDePolynomes_{\idp}\right)=\rg\left(V\AnneauDePolynomes_{\idp}\right)).
\end{equation}
\item Let $\mathcal{R}$ be a ring and $M$ be an $\mathcal{R}$-module. We denote by $l_{\mathcal{R}}(M)$ the length of $M$ (see p.~72 of~\cite{Eis} for the definition). In fact we shall use this definition only in the case $\mathcal{R}=\AnneauDePolynomes_{\idp}$ and $M=(\AnneauDePolynomes/I)_{\idp}$, where $I$ denotes an ideal of $\A$.
\item Let $I$ be a proper ideal of the ring $\AnneauDePolynomes$,
   \begin{equation} \label{definDePP_defin_m}
    m(I)=m(\eq(I)) \eqdef \sum_{\idp\in\Spec(\AnneauDePolynomes)\,\vline\,\rg(\idp)=\rg(I)}l_{\AnneauDePolynomes_{\idp}}((\AnneauDePolynomes/I)_{\idp}) \in \mnn^*.
   \end{equation}
\end{enumerate}
\end{definition}
}

The quantity $m(I)$ is the number of primary components of $I$ counted with their length as a multiplicity.

\begin{lemma} \label{LemmeCor14NumberW} 
Let $\idp\subset\AnneauDePolynomes$ be a prime bi-homogeneous ideal such that
$\V(\idp) \subset \mpp^1\times\mpp^n$ is projected onto $\mpp^1$. We recall the notations $\nu$, $V_i=V_i(\idp)$ and $\rho_i$ introduced in
Definition~\ref{V_irho_i} and $I(V_i,\idp)=(V_i\A_{\idp})\cap\A$ introduced in Definition~\ref{def_I}. One has the following upper bound for $m(I(V_i,\idp))$:
\begin{equation} \label{majorationm}
m(I(V_i,\idp)) \leq \nu(n+1)!\rho_i^{n+1}.
\end{equation}
\end{lemma}

The proof of this lemma needs quite many calculations. 
For technical reasons we have to distinguish two cases according to the value of $\nu_1$ in the condition~(\ref{degphiQleqdegQ}). 
More precisely, it is important in our proof to treat separately the cases where the constant $\nu_1$ is zero or non-zero.

We give a proof of Lemma~\ref{LemmeCor14NumberW} for the case $\nu_1=0$ in Section~\ref{section_proof_LemmeCor14NumberW}. This case is important in the treatment of differential operators (in particular for our improvement of Nesterenko's result). The proof in the case $\nu_1\ne 0$ is similar, it can be found in~\cite{EZ}, subsection~2.2.2.

\subsection{Definitions and properties related to multiprojective diophantine geometry} \label{definitions_multiprojective_dg}

We shall operate in the further text with such notions as \emph{(bi-)degree} and \emph{height} of a variety. Reader is invited to consult these notions in chapters~5 and~7 of~\cite{NP} or chapter~1 of~\cite{EZ}. We cite here just several properties that we shall use later.
\begin{enumerate}
  \item In case of hypersurface, i.e. if variery $V\subset\mpp^1_{\kk}\times\mpp^n_{\kk}$ is just a zero locus of a bi-homogeneous polynomial $P\in\A$, bi-degree of $V$ is a couple of integers $\left(\deg_{\ul{X}'}P,\deg_{\ul{X}}P\right)$. In this case it is common to denote also $\deg_{1,n-1}V:=\deg_{\ul{X}}P$  and $\deg_{0,n}V:=\deg_{\ul{X}'}P$. In general, bi-degree of variety $V$ is a couple of integers denoted often as $\left(\deg_{0,\dim(V)}V,\deg_{1,\dim(V)-1}V\right)$. This notation is explained in Chapter~5 of~\cite{NP}.
  \item If $V=V_1\cup\dots\cup V_r$ is a decomposition of $V$ in a union of irreducible components, we have
  \begin{equation*}
    \dim_{i,n-i}V=\sum_{j=1}^r\dim_{i,n-i}V_j,\quad i=0,1.
  \end{equation*}
  \item For any irreducible variety $V\subset\mpp^1\times\mpp^n$ and any hypersurface $Z\subset\mpp^1\times\mpp^n$ of bi-degree $(a,b)$, such that zero locus of $V$ is not a subset of $Z$ (i.e. $V$ and $Z$ intersect properly), there exists a variety $W$ such that its zero locus coincides with intersection of zero loci of $V$ and $Z$ (hence $\dim W=\dim V-1$), and $W$ satisfies
      \begin{eqnarray*}
        \dd_{(1,\dim(V)-2)}(W) &=& b\cdot\dd_{(1,\dim(V)-1)}(V), \\
        \dd_{(0,\dim(V)-1)}(W) &=& a\cdot\dd_{(1,\dim(V)-1)}(V) + b\cdot\dd_{(0,\dim(V))}(V).
      \end{eqnarray*}
      We shall denote such a variety $W$ as $V\cap Z$.
  \item Let $W\subset\mpp^n_{\kk(z)}$ be a subvariety. We can replace $(1:z)$ by $(X_0':X_1')$ transforming $W$ into a subvariety $\tilde{W}\subset\mpp^1_{\kk}\times\mpp^n_{\kk}$. In this case \emph{height} of $W$ can be expressed as $h(W)=\dd_{(0,\dim(\tilde{W}))}(\tilde{W})$ and \emph{degree} of $W$ is $\deg(W)=\dd_{(1,\dim(\tilde{W})-1)}(\tilde{W})$.
  \item To any bi-homogeneous ideal $I\subset\A$ (resp. any homogeneous ideal $J\subset\kk[z][X_0,\dots,X_n]$) we can associate a bi-projective (resp. projective) variety \V(I), thus defining $\deg_{i,n+1-\rk(I)-i}I$, $i=0,1$ (resp. $\deg(I)$ and $h(I)$).
\end{enumerate}

We shall regularly use the valuation $\ordz$ on the ring $\kk[[z]]$ of formal power series. This valuation induces the notions $\Ord({x},V)$ and $\ord({x},V)$, both measuring how far a point $x$ in a (multi-)projective space is from a variety $V$ emerged to the same space. Sometimes it is more convenient to use the quantities $\Dist({x},V):=\exp(-\Ord({x},V))$ and $\dist({x},V):=\exp(-\ord({x},V))$. Precise definitions could be found in~\cite{NP}, chapter~7, \S~4 and~\cite{EZ}, chapter~1, \S~3.

In order to make this article more self-containing we introduce briefly these notions (however for $\Ord$ we are forced to give just a short illustration, because technical details are quite complicated in this case).
\begin{definition} \label{defin_ord_xy} \begin{enumerate}
  \item If $\ul{x}=(x_0,\dots,x_n)\in\kk((z))^{n+1}$, we define $\ordz\ul{x}=\min_{i=0,\dots,n}\ordz x_i$.
  \item Let $x,y\in\mpp_{\kk((z))}^n$ be two points and $\ul{x}$ and \ul{y} be systems of projective coordinates respectively for $x$ and $y$. We define $x\wedge y$ to be a vector with $n(n-1)/2$ coordinates $\left(x_iy_j-x_jy_i\right)_{1\leq i<j\leq n}$ (the order of coordinates $x_iy_j-x_jy_i$ of this vector is not important for our purposes). Finally, we define
      \begin{equation} \label{def_ord_xy}
        \ordz(x,y):=\ordz(\ul{x}\wedge\ul{y})-\ordz \ul{x} - \ordz\ul{y}.
      \end{equation}
      One readily verifies that the r.h.s. in~(\ref{def_ord_xy}) does not depend on the choice of systems of projective coordinates for $x$ and $y$.
  \item Let $x,y\in\mpp^1_{\kk((z))}\times\mpp^n_{\kk((z))}$ and $\pi_1$ (resp. $\pi_n$) be a canonical projection of $\mpp^1_{\kk((z))}\times\mpp^n_{\kk((z))}$ to $\mpp^1_{\kk((z))}$ (resp. $\mpp^n_{\kk((z))}$). We define
       \begin{equation} \label{def_ord_xy_biprojectif}
        \ordz(x,y):=\min_{i=1,n}\ordz\left(\pi_i(x),\pi_i(y)\right).
      \end{equation}
  \item Let $V\subset\mpp^1\times\mpp^n$ (or $V\subset\mpp^n$) be a variety. We define
      \begin{equation} \label{def_ord_Vx}
        \ordz(x,V):=\max_{y\in V}\ordz\left(x,y\right).
      \end{equation}
  \item Sometimes we shall write simply $\ord(x,y)$, $\ord(x,V)$ and so on instead of $\ordz(x,y)$, $\ordz(x,V)$... This will not create any ambiguity because we shall be interested in only one valuation $\ordz$, so all the derived constructions, such as $\ord(x,y)$ and $\ord(x,V)$, will refer always to this valuation.
\end{enumerate}
\end{definition}
We proceed to introducing $\Ord(x,V)$. In this case however we give not a rigorous definition but rather an intuitive illustration (it concerns point~(\ref{defin_ordOrd_Ord2}) below). We refer the reader to~\cite{NP}, chapter~7, \S~4 and~\cite{EZ}, chapter~1, \S~3 for a rigorous (and strictly algebraic) definition of $\Ord(x,V)$.
\begin{enumerate}
  \item \label{defin_ordOrd_Ord1} First of all, if $V$ is 0-dimensional over $\ol{\kk((z))}$, it can be represented as a union of $r$ points $y_1,\dots,y_r$ (in fact, $r=\deg(V)$) and we define $\Ord(x,V):=\sum_{i=1}^r\Ord(x,y_i)$. In particular, if $V$ contains just one point over $\ol{\kk((z))}$, we have $\Ord(x,y)=\ord(x,y)$.
  \item \label{defin_ordOrd_Ord2} Let $V$ be a variety of dimension $d\geq 0$. We consider an intersection of $V$ with $d$ hyperplanes $U_1,\dots,U_d$ in general position, it is a 0-dimensional variety $V\cap U_1\cap\dots\cap U_d$. In particular, the quantity $\Ord(x,V\cap U_1\cap\dots\cap U_d)$ is well-defined as a function of $(U_1,\dots,U_d)$. We denote $\Ord(x,V)$ the supremum of this function for all the sets of hyperplanes $U_1,\dots,U_d$ in general position.
  \item More generally, we can replace in point~(\ref{defin_ordOrd_Ord2}) a set of $d$ hyperplanes in general position by a set of $d$ hypersurfaces, of degrees $(\delta_1,\dots,\delta_d)\in\mnn^n$, in general position. In this case we obtain a notion usually denoted as $\Ord_{(\delta_1,\dots,\delta_d)}(x,V)$.
\end{enumerate}

We mention here two special properties that we shall use later. Proofs are easy and can be found in~\cite{EZ}, Chapter~1.

\begin{lemma} \label{Representants} Let $\b{x},\b{y} \in \mpp^n_{\overline{\kk((\b{z}))}}$ be two points in the projective space.

a) Let $\ull{x}$ be a system of projective coordinates of $\b{x}$ and $\ull{y}$ be a system of projective coordinates of $\b{y}$ satisfying
\begin{equation*}
 \ordz \ull{x} = \ordz \ull{y}.
\end{equation*}
Then
\begin{equation} \label{LemmeRepresentantsA}
 \Ordz(\b{x},\b{y}) \geq \ordz(\ull{x} - \ull{y}) - \ordz\ull{x}.
\end{equation}

b) Suppose $\Ordz(\b{x},\b{y})>0$, if we fix for $\b{y}$ a system system of projective coordinates $\ull{y}$ in $\overline{\kk((\b{z}))}^{n+1}$, then there is a system of projective coordinates $\ull{x} \in \overline{\kk((\b{z}))}^{n+1}$ of $\b{x}$ satisfying
\begin{equation}
 \begin{split}
  &\alpha) \quad \ordz\ull{x}=\ordz\ull{y},\\
  &\beta) \quad\ordz(\ull{x} - \ull{y}) - \ordz(\ull{y}) = \Ordz(\b{x},\b{y})
\end{split}
\end{equation}
\end{lemma}
\begin{proof}
See Lemma~1.22 of~\cite{EZ}.
\end{proof}

\begin{lemma}[Liouville's inequality]
Let $Q\in\kk(\b{z})$ and $\b{Z}$ be a cycle of dimension 0 defined over $\kk(\b{z})$. Then
\begin{equation}\label{iet_main}
    \deg(\b{Q})h(\b{Z})+h(\b{Q})\deg(\b{Z})\geq\left|\sum_{\b{\beta}\in\b{Z}}\ord_{v_0}\left(\b{Q}(\ull{\beta})\right)\right|,
\end{equation}
\end{lemma}
\begin{proof}
See inequality~(1.20) at the end of section~1.2.2 of~\cite{EZ}.
\end{proof}

In the following definition we associate to each bi-projective ideal $I$ a couple of integers, $\left(\delta_0(I),\delta_1(I)\right)$. This quantity plays an important role in our article (though it seems to be a little bit too complicated at first glance), so we make a short intuitive comment in Remark~\ref{rem_def_delta} just after the definition.
\begin{definition} \label{def_delta}
a) Let $I \subset \AnneauDePolynomes$ be a bi-homogeneous ideal, we choose a bi-homogeneous polynomial $P \in I \setminus \{0\}$
that minimizes the quantity
\begin{equation} \label{def_delta_condition_minimum_crossproduct}
\mu\dd_{(0,n-\rg I+1)}I\deg_{\ul{X}}P + \nu_0\dd_{(1,n-\rg I)}I \deg_{\ul{X}'}P + \nu_1\dd_{(1,n-\rg I)}I \deg_{\ul{X}}P.
\end{equation}
If there exists more than one bi-homogeneous polynomial minimising~(\ref{def_delta_condition_minimum_crossproduct})
we choose a polynomial with $\deg_{\ul{X}'}P$ minimal.
We introduce notation $\delta_0(I) \eqdef \deg_{\ul{X}'}P$ and $\delta_1(I) \eqdef \deg_{\ul{X}}P$.

b) For all cycle $Z$ (defined over $\kk$) in $\mpp^1\times\mpp^n$
we define $\delta_i(Z) \eqdef \delta_i(\I(Z))$, $i=0,1$.
\end{definition}
\begin{remark} \label{rem_def_delta}
To each (homogeneous, say) ideal $I$ in the ring of polynomials we can associate an integer $\delta$, the minimum of degree of polynomials belonging to $I$:
\begin{equation*}
    \delta:=\min_{P\in I}\deg(P).
\end{equation*}
In bi-homogeneous case degree is replaced by couple of integers $\left(\delta_0(I),\delta_1(I)\right)$. For some reasons (basically, the use of Bezout's theorem, e.g. Theorem~4.11 of chapter~3 of~\cite{NP}) the use of a simple expression such as
\begin{equation}\label{def_delta_naive}
    \min_{P\in I}\left(\deg_{X'}P+\deg_{X}P\right)
\end{equation}
does not lead to good results. The correct generalization in our case needs introducing weights for $\deg_{X'}$ and $\deg_{X}$ depending on the ideal $I$ as well as constants $\mu$, $\nu_0$, $\nu_1$ characterizing the transformation $\phi$ (see the property~(\ref{degphiQleqdegQ})). This "weighted" version of~(\ref{def_delta_naive}) is realized in~(\ref{def_delta_condition_minimum_crossproduct}). Although the obtained notion seems to be quite sophisticated, it allows to perform proofs presented in this article.
\end{remark}
Here is the first property of the quantity $\left(\delta_0(I),\delta_1(I)\right)$.
\begin{lemma}\label{deltainfty} Fix a point
$$
\ul{f}=(1:\b{z},1:f_1(\b{z}):...:f_n(\b{z})) \in \mpp^1_{\kk[[\b{z}]]}\times\mpp^n_{\kk[[\b{z}]]}
$$
and consider a sequence of cycles $Z_i \subset \mpp^1\times\mpp^n$ defined over $\kk$ and such that $\ul{f} \not\in Z_i$ for $i=1,2,3,\dots$. If $\ord(\ul{f},Z_i)$ tends to $+\infty$ (as $i\rightarrow\infty$), then $\max(\delta_0(Z_i,\ul{f}),\delta_1(Z_i,\ul{f}))$ also tends to the infinity (as $i\rightarrow\infty$).
\end{lemma}
The proof of Lemma~\ref{deltainfty} is easy and can be found in~\cite{EZ} (Lemme~1.23). We shall need only its weak corollary:

\begin{corollary}\label{cor1_deltainfty} There exist a constant $C_{sg}$ which depends only on $\ul{f}=(1:\b{z},1:f_1(\b{z}):...:f_n(\b{z})) \in \mpp^1_{\kk[[\b{z}]]}\times\mpp^n_{\kk[[\b{z}]]}$ such that if a cycle $Z \subset \mpp^1\times\mpp^n$ (defined over $\kk$) does not contain $\ul{f}$ and satisfies $\ord_{\ul{f}}Z \geq C_{sg}$, then either $\delta_0(Z,\ullt{f})\geq 2n!+1$ or $\delta_1(Z,\ullt{f})\geq 4n$.
\end{corollary}

We introduce now a key notion of forthcoming proofs. As before, we discuss the notions introduced here in Remark~\ref{rem_V_irho_i} just after the definition
\begin{definition} \label{V_irho_i} We define
\begin{equation}\label{def_nu}
    \nu\eqdef\begin{cases}1 &\text{ if } \nu_1=0,\\
                    2^{n+2}\max\left(1,\frac{4\nu_0}{\nu_1}\right)^{n+1} &\text{ if } \nu_1\ne 0.
             \end{cases}
\end{equation}
The sequence of numbers $\rho_i$ is defined recursively. We put $\rho_0=0$, $\rho_1=1$ and $\rho_{i+1} = \nu(n+1)!\rho_{i}^{n+2}\max\left(\mu,\nu_0\right)^{\nu(n+1)!\rho_{i}^{n+1}}$
for $i=1,...,n+1$. The constants $\mu$, $\nu_0$ et $\nu_1$ in this definition are the same as in~(\ref{degphiQleqdegQ}).

Let $Z$ be an algebraic cycle (defined over $\kk$) in the space $\mpp^1\times\mpp^n$.
We denote by $V_i$, or more precisely by $V_i(Z)$, the vectorial space (over $\kk$) generated by the
polynomials (in $\kk[X_0':X_1',X_0:...:X_n]$) vanishing over the cycle $Z$ and of the degree in $\ul{X}'$ at most $\rho_i\left(\delta_0(Z)+\frac{\nu_1}{\max(\mu,\nu_0)}\delta_1(Z)\right)$ and of the degree in $\ul{X}$ at most $\rho_i\delta_1(Z)$
(recall that $\delta_0(Z)$ and $\delta_1(Z)$ are introduced in Definition~\ref{def_delta}).

If $I$ is a proper bi-homogeneous ideal of $\AnneauDePolynomes$ we also use the notation
\begin{equation*}
    V_i(I)=V_i(\V(I)),
\end{equation*}
where $\V(I)$ is the cycle of $\mpp^1\times\mpp^n$ defined by $I$.
\end{definition}
\begin{remark} \label{rem_V_irho_i}
Basically, $V_i$ represents a set of polynomials of bi-degree only a constant ($\rho_i$) times bigger than a "minimal" bi-degree (in terms explained in Definition~\ref{def_delta} and Remark~\ref{rem_def_delta}). Note that constants $\rho_i$ are absolute (they depend only on initial data: number of functions $n$ and constants $\mu$, $\nu_0$ and $\nu_1$ controlling the growth of degree of polynomial under the action of $\phi$).

\end{remark}


\begin{definition} \label{def_i0} We associate to every variety $W$ a number $i_0(W)$. We define it to be the biggest positive integer such that
$\rg\left(I_0\left(V_i(W),\I(W)\right)\right)\geq i$ for all $1\leq i\leq i_0(W)$.
\end{definition}


\begin{remark} \label{rem_i0}
Definitions~\ref{def_delta} and~\ref{V_irho_i} been given, one verifies the inequality $$\rg\left(I_0\left(V_1(W),\I(W)\right)\right)\geq 1.$$ So, the index $i_0(W)\geq 1$ is well defined
for all varieties $W$. On the other hand the rank of any homogeneous ideal in $\AnneauDePolynomes$ can not exceed $n+1$, thus $i_0(W)\leq n+1$ for all variety $W$.
\end{remark}

\section{Proof of Lemma~\ref{LemmeCor14NumberW}} \label{section_proof_LemmeCor14NumberW}

In this section we prove Lemma~\ref{LemmeCor14NumberW} for the case $\nu_1=0$. We shall deduce it as a corollary of the following technical lemma (Lemma~\ref{LemmeProp14_1}).
The proof in the case $\nu_1\ne 0$ is similar and can be found in~\cite{EZ}, subsection~2.2.2.

\begin{lemma} \label{LemmeProp14_1} Let $I$ be a bi-homogeneous ideal of $\AnneauDePolynomes$, $I\ne\AnneauDePolynomes$ and $\delta_0=\delta_0(I)$, $\delta_1=\delta_1(I)$ are the quantities introduced in Definition~\ref{def_delta}.

Let $W \subsetneq \mpp^1\times\mpp^n$ be an irreducible variety
 projected onto $\mpp^1$ and such that
\begin{equation} \label{LemmeProp14_hypothese1_1}
\delta_1\dd_{(0,\dim W)}W + \delta_0\dim W\dd_{(1,\dim W-1)}W <
\frac{\delta_0\delta_1^{\codim(W)}}{n!}.
\end{equation}
Then there is a polynomial $Q \in \I(W) \setminus I$ satisfying two inequalities:
\begin{equation} \label{LemmeProp14_petitPolynome_1}
\begin{aligned}
\deg_{\ul{X}'}Q&\leq\delta_0-1,\\
\deg_{\ul{X}}Q&\leq\delta_1-1.
\end{aligned}
\end{equation}
\end{lemma}

\begin{proof} Note that in the case $\delta_0=0$ or $\delta_1=0$
the condition~(\ref{LemmeProp14_hypothese1_1}) can not be satisfied: the r.h.s. is zero and can not exceed strictly the l.h.s.
So it is enough to consider the case $\delta_0 \geq 1$ and $\delta_1 \geq 1$. Since $W$ is projected onto $\mpp^1$, we have
\begin{equation} \label{LemmeProp14_degW_geq_1_1}
\dd_{(1,\dim W-1)} W \geq 1,
\end{equation}
so the l.h.s. of~(\ref{LemmeProp14_hypothese1_1}) is at least $\dim W\geq 1$. 
Hence we can suppose
$$
\delta_0\delta_1^{\codim(W)} \geq n!
$$

If a polynomial $Q_0$ of bi-degree $(\delta_0-1,\delta_1-1)$ vanishes on $W$, it satisfies all the announced properties of Lemma~\ref{LemmeProp14_1}, there is nothing more to prove in this case. Indeed, $Q_0$ does not belong to $I$ because
\begin{multline} \label{LemmeProp14_ie1_1}
\mu(\delta_1-1)\dd_{(0,\dim I)} I + \nu_0(\delta_0-1)\dd_{(1,\dim I-1)}I + \nu_1(\delta_1-1)\dd_{(1,\dim I-1)}I\\<\mu\delta_1\dd_{(0,\dim I)} I + \nu_0\delta_0\dd_{(1,\dim I-1)}I
 + \nu_1\delta_1\dd_{(1,\dim I-1)}I
\end{multline}
and the r.h.s. of~(\ref{LemmeProp14_ie1_1}) minimizes the expression
\begin{equation*}
    \mu\deg_{\ul{X}}(Q)\dd_{(0,\dim I)} I + \nu_0\deg_{\ul{X}'}(Q)\dd_{(1,\dim I-1)}I + \nu_1\deg_{\ul{X}}(Q)\dd_{(1,\dim I-1)}I
\end{equation*}
for all the bi-homogeneous polynomials $Q$ from $I\setminus\{0\}$ (see Definition~\ref{def_delta}).

It remains to eliminate the case when no polynomial of bi-degree $(\delta_0-1,\delta_1-1)$ vanishes on $W$.
Under this hypothesis, we have the estimation
\begin{equation} \label{LemmeProp14_Hgminoration_1}
H_g(W,\delta_0-1,\delta_1-1) \geq \delta_0\binom{\delta_1+n-1}{n} \geq \delta_0\frac{\delta_1^{n}}{n!}.
\end{equation}

In the same time, by Corollary~9 of~\cite{PP} we have
\begin{equation} \label{LemmeProp14_Hgmajoration_1}
\begin{aligned}
H_g&(W,\delta_0-1,\delta_1-1) \leq (\delta_1-1)^{\dim W}\dd_{(0,\dim W)}W \\&+ (\delta_0-1)(\delta_1-1)^{\dim W-1}\dim W\dd_{(1,\dim W-1)} W + \dim W.
\end{aligned}
\end{equation}
Since we have~(\ref{LemmeProp14_degW_geq_1_1}), we can weaken~(\ref{LemmeProp14_Hgmajoration_1}) to
\begin{equation} \label{LemmeProp14_Hgmajoration2_1}
\begin{aligned}
H_g(W,\delta_0-1,\delta_1-1) \leq &\delta_1^{\dim W}\dd_{(0,\dim W)}W \\&+ \delta_0\delta_1^{\dim W-1}\dim W\dd_{(1,\dim W-1)} W.
\end{aligned}
\end{equation}

Inequalities~(\ref{LemmeProp14_Hgminoration_1}) and~(\ref{LemmeProp14_Hgmajoration2_1}) put together contradict
assumption~(\ref{LemmeProp14_hypothese1_1}).
\end{proof}

\begin{corollary} \label{LemmeCor14_1} In the situation of lemma~\ref{LemmeProp14_1} if $I$ is a radical ideal and if an irreducible variety $W$ is projected
onto $\mpp^1$ and contains $\V(I)$,
then
\begin{equation} \label{minorationdegW_1}
\delta_1\dd_{(0,\dim W)}W + \delta_0\dim W\dd_{(1,\dim W-1)}W \geq
\frac{\delta_0\delta_1^{\codim(W)}}{n!},
\end{equation}
\end{corollary}
\begin{proof}
The proof is straightforward using Lemma~\ref{LemmeProp14_1}.
\end{proof}
\begin{proof}[Proof of Lemma~\ref{LemmeCor14NumberW}]
We let $r$ denote $\rg I_0(V_i,\idp)$. The ideal $I_0(V_i,\idp)\subset\idp$ is extended-contracted by localization at $\idp$ of an ideal generated by
polynomials of bi-degree $\left(\rho_i\delta_0(\idp),\rho_i\delta_1(\idp)\right)$ (see Definition~\ref{def_i0}). Hence it satisfies, according to B\'ezout's theorem,
\begin{multline} \label{majorationTordue0_1}
(n-r+1)(\rho_i\delta_0(\idp))(\rho_i\delta_1(\idp))^{n-r}\dd_{(1,n-r)}\left( I_0(V_i,\idp)\right)
\\+(\rho_i\delta_1(\idp))^{n-r+1}\dd_{(0,n-r+1)}\left( I_0(V_i,\idp)\right)\\ \leq(n+1)\rho_i^{n+1}\delta_0(\idp)(\delta_1(\idp))^{n}.
\end{multline}
By simplifying out the factor $\rho_i^{n-r+1}\delta_1(\idp)^{n-r}$ from the both sides of the latter inequality we obtain
\begin{multline} \label{majorationTordue_1}
(n-r+1)\delta_0(\idp)\dd_{(1,n-r)}\left( I_0(V_i,\idp)\right)
+\delta_1(\idp)\dd_{(0,n-r+1)}\left( I_0(V_i,\idp)\right)\\ \leq(n+1)\delta_0(\idp)(\rho_i\delta_1(\idp))^{r}.
\end{multline}
Let $W=\V(\idq)$, where $\idq$ is a prime ideal (not necessarily the minimal one) associated to $\V(I_0(V_i,\idp))$. By construction of $I_0(V_i,\idp)$
all its associated primes are contained in $\idp$, thus $\V(\idp) \subset W$. Moreover, since
$\V(\idp)$ is projected onto $\mpp^1$, we deduce that $W$ also is projected onto $\mpp^1$.
Applying Corollary~\ref{LemmeCor14_1} we obtain that $W$
satisfies~(\ref{minorationdegW_1}); in other words, every prime $\idq$ associated to $I_0(V_i,\idp)$ satisfies
\begin{equation} \label{minorationDegQ_1}
\delta_1(\idp)\dd_{(0,n+1-\rg \idq)}\idq + (n+1-\rg\idq)\delta_0(\idp)\dd_{(1,n-\rg\idq)}\idq\\ \geq
\frac{\delta_0(\idp)\delta_1(\idp)^{\rg\idq}}{n!}.
\end{equation}
But
\begin{equation} \label{degDecomposition_1}
\begin{aligned}
 \dd_{(1,n-r)}(I_0(V_i,\idp)) = \sum_{\substack{\idq \in \Spec\A,\\ \rg(\idq)=r}}\dd_{(1,n-r)}(\idq)l_{\AnneauDePolynomes_{\idq}}(\left(\AnneauDePolynomes/I_0(V_i,\idp)\right)_{\idq}),
\end{aligned}
\end{equation}
and
\begin{equation} \label{htDecomposition_1}
\begin{aligned}
 \dd_{(0,n-r+1)}(I_0(V_i,\idp)) = \sum_{\substack{\idq \in \Spec\A,\\ \rg(\idq)=r}}\dd_{(0,n-r+1)}(\idq)l_{\AnneauDePolynomes_{\idq}}(\left(\AnneauDePolynomes/I_0(V_i,\idp)\right)_{\idq}),
\end{aligned}
\end{equation}
Summing up~(\ref{degDecomposition_1}) with coefficient $(n-r+1)\delta_0(\idp)$
and~(\ref{htDecomposition_1}) with coefficient $\delta_1(\idp)$, we obtain
\begin{equation}\label{minorationDegHtRightParm_1}
\begin{aligned}
 &(n-r+1)\delta_0(\idp)\dd_{(1,n-r)}(I_0(V_i,\idp))+\delta_1(\idp)\dd_{(0,n-r+1)}(I_0(V_i,\idp))\\
 &=\sum_{\substack{\idq \in \Spec\A,\\ \rg(\idq)=r}}\left((n-r+1)\delta_0(\idp)\dd_{(1,n-r)}\idq + \delta_1(\idp)\dd_{(0,n-r+1)}\idq\right)\\
 &\times l_{\AnneauDePolynomes_{\idq}}(\left(\AnneauDePolynomes/I_0(V_i,\idp)\right)_{\idq}).
\end{aligned}
\end{equation}
Applying~(\ref{majorationTordue_1}) to the l.h.s. of~(\ref{minorationDegHtRightParm_1}) and~(\ref{minorationDegQ_1}) to the r.h.s. of~(\ref{minorationDegHtRightParm_1}) we obtain
\begin{equation} \label{majorationm_pre_1}
(n+1)\delta_0(\idp)(\rho_i\delta_1(\idp))^{r} \geq \frac{\delta_0(\idp)\delta_1(\idp)^{r}}{n!}m(I_0(V_i,\idp))
\end{equation}
Finally, we deduce~(\ref{majorationm}) from~(\ref{majorationm_pre_1}) with simplification and using the remark $r \leq n+1$.
\end{proof}

\section{Transference lemma of P.Philippon}

In this section we present (a particular case of) the transference lemma elaborated in~\cite{PP}.

In the sequel we shall use the constant $c_n$ defined by
\begin{equation} \label{def_cn}
c_n=2^{n+1}(n+2)^{(n+1)(n+3)}.
\end{equation}
Note that in terms of~\cite{PP} one has $c_n=c_{\mpp^1\times\mpp^n}$.

\begin{theorem}[Transference $(1,n)$-Projective Lemma]\label{LdT}
Let $\ull{f} \in \mpp^n_{\kk[[\b{z}]]}$ and $C$ be a real number satisfying 
\begin{equation} \label{LdT_Cestgrande}
C\geq\left(\min(\nu_0,\mu)^{n}\right)^{-1}.
\end{equation}
If a non-zero form $\b{P} \in \kk[\b{z}][X_0,X_1,...,X_n]$ satisfies
\begin{multline} \label{ordPplusqueb}
 \ordz (\b{P}(\ull{f})) - \deg\b{P}\cdot\ordz (\ull{f}) - h(\b{P})\\ > C\cdot n\cdot\left((\mu+\nu_0) \left(h(\b{P})+1\right)+\nu_1\deg\b{P}\right)\mu^{n-1}\left(\deg\b{P}+1\right)^{n},
\end{multline}
then there is an  irreducible cycle $\b{Z} \in \mpp_n\left(\overline{\kk(\b{z})}\right)$ defined over $\kk(\b{z})$, of
dimension~0, contained in the zero locus of $\b{P}$, satisfying
\begin{multline} \label{LdTdegZ}
 \nu_0\deg\b{Z}\cdot h(\b{P})+\nu_1\deg\b{Z}\cdot\deg\b{P}+\mu\cdot h(\b{Z})\cdot\deg\b{P}\\ \leq c_n C(n+1) \mu^{n} \left(\nu_0\left(h(\b{P})+1\right)+\nu_1\deg\b{P}\right)\left(\deg\b{P}+1\right)^{n},
\end{multline}
and
\begin{equation} \label{LdTordZ0}
\begin{split}
\sum_{\alpha \in \b{Z}} \Ord_{\ull{f}}(\alpha) > &c_n^{-1}\left[(c_nC)^{\frac{1}{n}}\left(\nu_0(h(\b{P})+1)+\nu_1\deg\b{P}\right)\right]\deg(\b{Z})\\
&+c_n^{-1}\left[(c_nC)^{\frac{1}{n}}\mu(\deg\b{P}+1)\right]h(\b{Z}).
\end{split}
\end{equation}
In particular, (\ref{LdT_Cestgrande}) and~(\ref{LdTordZ0}) imply
\begin{equation} \label{LdTordZ}
\begin{split}
\sum_{\ul{\alpha} \in \b{Z}} \Ord_{\ull{f}}(\ul{\alpha}) > C^{\frac{1}{n}}c_n^{-\frac{n-1}{n}}\Big( \nu_0\deg(\b{Z})h(\b{P})&+\nu_1\deg(\b{Z})\deg\b{P}\\&+ \mu\cdot h(\b{Z})\deg\b{P} \Big).
\end{split}
\end{equation}
\end{theorem}
\begin{proof}
We apply Corollary~11 of~\cite{PP} to $\tilde{X}_0=\V(\b{P})$ and $\tilde{\Phi}=\ullt{f}$ with
\begin{eqnarray*}
  \eta &=& \left[(c_nC)^{\frac{1}{n}}\left(\nu_0(h(\b{P})+1)+\nu_1\deg\b{P}\right)\right], \\
  \delta &=& \left[(c_nC)^{\frac{1}{n}}\mu(\deg\b{P}+1)\right].
\end{eqnarray*}
Inequality~(\ref{LdT_Cestgrande}) implies
\begin{equation} \label{LdT_preuve_majoration_hdeg}
\begin{aligned}
   h(\b{P}) &\leq&\eta \\
   \deg\b{P} &\leq&\delta,
\end{aligned}
\end{equation}
thus $\tilde{X}_0$ id defined by a form of multidegree $\leq(\eta,\delta)$ with
$$
(\eta,\delta)\geq(c_n^{1/n},c_n^{1/n}).
$$
In our case $c_n=c_{\mpp^1\times\mpp^n}$ and
\begin{equation} \label{LdT_preuve_degX0etadelta}
\deg\left(\tilde{X}_0;(\eta,\delta)\right)=h(\b{P})\cdot\delta^n+n\cdot\deg\b{P}\cdot\eta\delta^{n-1}.
\end{equation}
The condition
$$
\Ord_{\ull{f}}\tilde{X}_0\geq c_n^{-1}\deg\left(\tilde{X}_0;(\eta,\delta)\right)
$$
is assured by~(\ref{ordPplusqueb}), because
$$
\begin{aligned}
\Ord_{\ullt{f}}\tilde{X}_0&=\Ord_{\ull{f}}\b{P}=\ordz (\b{P}(\ull{f})) - \deg_{\ul{X}}\cdot\b{P} \ordz (\ull{f}) - h(\b{P})\\
&>C\cdot n\cdot\left((\mu+\nu_0) \left(h(\b{P})+1\right)+\nu_1\deg_{\ul{X}}\b{P}\right)\mu^{n-1}\left(\deg_{\ul{X}}\b{P}+1\right)^{n}\\
&>c_n^{-1}\cdot n\cdot(c_nC)\Bigg(\left(\nu_0(h(\b{P})+1)+\nu_1\deg_{\ul{X}}\b{P}\right)\mu^{n-1}\left(\deg_{\ul{X}}\b{P}+1\right)^{n}\\
&\quad + (h(\b{P})+1)\mu^n\left(\deg_{\ul{X}}\b{P}+1\right)^n\Bigg)\\
&> c_n^{-1}\Big(n\cdot\eta\delta^{n-1}\deg_{\ul{X}}\b{P}+\delta^nh(\b{P})\Big)\\
&>c_n^{-1}\deg\left(\tilde{X}_0;(\eta,\delta)\right).
\end{aligned}
$$
The conclusion of Corollary~11 of~\cite{NP} gives us exactly the conclusion of the theorem. Indeed, this corollary provides us a cycle $\b{Z}\subset\tilde{X}_0(\ol{\kk(\b{z})})$ defined over $\kk(\b{z})$ and of dimension 0 such that
\begin{eqnarray}
  \delta\cdot h(\b{Z})+\eta\deg\b{Z} &\leq& \deg\left(\tilde{X}_0,(\eta,\delta)\right) \label{LdT_preuve_cor_conclusion_deg},\\
  \sum_{\ul{\alpha}\in\b{Z}} \Ord_{\ull{f}}(\ul{\alpha}) &>& c_n^{-1}\left(\eta\deg\b{Z} + \delta\cdot h(\b{Z})\right) \label{LdT_preuve_cor_conclusion_Ord}.
\end{eqnarray}
Inequality~(\ref{LdT_preuve_cor_conclusion_deg}) (together with~(\ref{LdT_preuve_majoration_hdeg}) and~(\ref{LdT_preuve_degX0etadelta})) gives us inequality~(\ref{LdTdegZ}), and~(\ref{LdT_preuve_cor_conclusion_Ord}) provides us~(\ref{LdTordZ0}).
\end{proof}
\begin{definition} \label{defZP} Let $C$ be a real number satisfying~(\ref{LdT_Cestgrande}). We associate to each homogeneous polynomial $\b{P} \in \kk[\b{z}][X_1,...,X_n]$ and a real constant $C>0$ satisfying~(\ref{ordPplusqueb}) 
an irreducible 0-dimensional cycle $\b{Z}_{C}(\b{P})$ defined over $\kk(\b{z})$, contained in the zero locus of $P$ and satisfying
inequalities~(\ref{LdTdegZ}) and~(\ref{LdTordZ0}). In view of the transference lemma there exists at least one cycle satisfying all these conditions (provided polynomial $P$ and constant $C$ satisfy~(\ref{ordPplusqueb})). If there exists more than one such cycle, we choose one of them and fix this choice.
\end{definition}

\begin{remark} \label{genZP} Considering $(1:\b{z})$ as coordinates of a point in $\mpp^1_{\overline{\kk(\b{z})}}$ we can consider the cycle $\b{Z}$ as a $1$-dimensional cycle in $\mpp^1_{\kk}\times\mpp^n_{\kk}$ (defined over $\kk$). In this case we denote this cycle by $\Z_C(P)$.

\smallskip

We associate to a bi-homogeneous polynomial $P(X_0',X_1',X_0,X_1,...,X_n)\in\AnneauDePolynomes$ satisfying
\begin{equation} \label{ordPplusqueb2}
\frac{\ordz (P(1,\b{z},\ull{f}) - (\deg_{\ul{X}}P) \ordz (\ull{f}) - \deg_{\ul{X}'}P}{n((\nu_0+\mu)(\deg_{\ul{X}'}P+1)+\nu_1\deg_{\ul{X}}P)\mu^{n-1}(\deg_{\ul{X}}P+1)^{n}} > C,
\end{equation}
the homogeneous polynomial
$$
\b{P}(X_0,X_1,...,X_n)=P(1,\b{z},X_0,X_1,...,X_n)
$$
(satisfying in this case~(\ref{ordPplusqueb})). We have already defined the cycles $\b{Z}_C(\b{P})$ et $\Z_{C}(\b{P})$ for the latter polynomial. By this procedure we associate equally the cycles $\b{Z}_C(P)$ and $\Z_{C}(P)$ to every bi-homogeneous polynomial $P \in \AnneauDePolynomes$ satisfying~(\ref{ordPplusqueb2}).
\end{remark}

\begin{remark} \label{rem_LdT_Z_nonisotrivial} Note that combining~(\ref{ordPplusqueb2}) (for $C$ large enough) with the transference lemma (Theorem~\ref{LdT}, (\ref{LdTordZ0})) one shows that the cycle $\Z_{C}(P)$ is not an isotrivial one (and thus $\b{Z}_C(P)$ is not defined over~$\kk$). Indeed, it is easy to find a lower bound for $C$ in order to assure that $\Z_{C}(P)$ is not isotrivial. Each point defined over $\kk$ contributes at most $\Ordz\left(\ullt{f}\wedge\ull{f}(0)\right)$ to $\Ord_{\ull{f}}\Z_{C}(P)$, so for an isotrivial cycle $Z$ one has
$$
\Ord_{\ull{f}}Z\leq\Ordz\left(\ull{f}(\b{z})\wedge\ull{f}(0)\right)\deg Z.
$$
Thus
$$
C>\Ciso:=\left(\frac{c_n\Ordz\left(\ull{f}\wedge\ull{f}(0)\right)+1}{\min(\nu_0,\mu)}\right)^n
$$
implies that $\Z_{C}(P)$ is not isotrivial.
\end{remark}

The following theorem plays an important role in the proof of our principal result, Theorem~\ref{LMGP}. More precisely, without Theorem~\ref{dist_alpha} we could prove (with essentially the same arguments) Theorem~\ref{LMGP} with the assumption~(\ref{RelMinN2}) replaced by a stronger hypothesis:
\begin{equation*}
    \ord_{\ull{f}}(\idq) < K_0.
\end{equation*}
Note that corresponding corollary for differential operators in the case of characteristic~0 coincides with Nesterenko's conditional multiplicity lemma (see Theorem~1.1 of Chapter~10 \cite{NP}, or Theorem~\ref{theoNesterenko_classique} below).

\begin{theorem} \label{dist_alpha} Let $\ull{f}=(1:\b{f}_1:...:\b{f}_n) \in
\mpp^n_{\kk[[\b{z}]]}$ be a point with all the coordinates $\b{f}_1$, $\b{f}_2$,...,$\b{f}_n$ algebraically independent
over $\kk(\b{z})$. Let $\b{P} \in
\kk[\b{z}][X_0,...,X_n]$ be a non-zero homogeneous polynomial satisfying~(\ref{ordPplusqueb}) with
\begin{equation} \label{dist_alpha_Cestgrande}
C \geq \max\left(\left(3/\min(\nu_0,\mu)\right)^n(n!)^nc_n^{n-1},\left(\frac{c_nC_{sg}+1}{\min(\nu_0,\mu)}\right)^n,\Ciso\right)
\end{equation}
(where $C_{sg}$ is the constant from Corollary~\ref{cor1_deltainfty}).
Let $\b{Z}=\b{Z}_C(\b{P})$ and let $\b{P}_0 \in \kk[\b{z}][X_0,...,X_n]$ be a homogeneous non-zero polynomial in $\ul{X}$, vanishing on $\b{Z}$ and realizing the minimum of the quantity
\begin{equation} \label{dist_alpha_q}
 \nu_0\cdot h(\b{Z})h(\b{Q})+\nu_1\cdot\deg\b{Z}\cdot h(\b{Q})+\mu\cdot h(\b{Z})\deg_{\ul{X}}\b{Q}
\end{equation}
over all the homogeneous polynomials $\b{Q}\in\kk[\b{z}][X_0,...,X_n]$ different from 0 and vanishing on $\b{Z}$.

We introduce the notation  $\delta_0=h(\b{P}_0)$ and $\delta_1=\deg_{\ul{X}}\b{P}_0$.

Then there exists a point $\ull{\alpha} \in \b{Z}$ satisfying
\begin{equation} \label{Cdirect}
\Ord(\ull{f},\ull{\alpha})
> \tilde{C} (\delta_0+1)(\delta_1+1)^n,
\end{equation}
where  $\tilde{C}=C^{\frac{1}{n}}\min(\nu_0,\mu)\left(3\cdot n!\cdot c_n^{\frac{n-1}{n}}\right)^{-1}$.
\end{theorem}
\begin{proof}
We claim that there exists a point $\ull{\alpha}_1 \in \b{Z}_C(\b{P})$ satisfying
\begin{equation} \label{dist_alpha_Ord_assezgrand}
\Ord(\ull{\alpha}_1,\ull{f}) \geq C_{sg}.
\end{equation}
Indeed, the lower bound~(\ref{LdTordZ0}) implies that there exists a point $\ull{\alpha}_1 \in \b{Z}_C(\b{P})$ satisfying $\Ord(\ull{\alpha}_1,\ull{f}) \geq c_n^{-1}\left[(c_nC)^{\frac{1}{n}}\nu_0\right]$ and we deduce~(\ref{dist_alpha_Ord_assezgrand}) from~(\ref{dist_alpha_Cestgrande}).

In view of Corollary~\ref{cor1_deltainfty} we can thus suppose
\begin{equation} \label{dist_alpha_soit0_soit1}
\delta_0 \geq 2\cdot n!+1 \mbox{ or } \delta_1 \geq 4n.
\end{equation}

Let $(a,b)\in\mnn^2$. By linear algebra one can construct a bi-homogeneous
polynomial $Q_{(a,b)}=Q_{(a,b)}(X_0',X_1',X_0,X_1,...,X_n)\in\AnneauDePolynomes
\setminus \{0\}$ of bi-degree $(a,b)$ and of vanishing order at
$\ullt{f}=(1,\b{z},\ull{f})$ satisfying
\begin{equation}\label{ie_Qab}
    \ordz Q_{(a,b)}(\ullt{f}) \geq \lfloor\frac{1}{n!}(a+1)(b+1)^n\rfloor.
\end{equation}
Let
\begin{equation} \label{dist_alpha_choix_ab}
(a,b)=\left\{\begin{aligned}&(\delta_0-1,\delta_1), \; \mbox{ if }\delta_0 \geq 2\cdot n!+1,\\
                     &(\delta_0,\delta_1-1), \; \mbox{ otherwise, i.e. } \delta_1 \geq 4n \mbox{ in view of~(\ref{dist_alpha_soit0_soit1}) }.\end{aligned}\right.
\end{equation}
We claim that for this choice of $(a,b)$ the following inequality holds
\begin{equation} \label{dist_alpha_est_ab}
\ordz Q_{(a,b)}(\ullt{f}) > \frac{1}{2\cdot n!}(\delta_0+1)(\delta_1+1)^n.
\end{equation}
In view of~(\ref{dist_alpha_choix_ab}), it suffices to consider two cases:
\begin{list}{\alph{tmpabcd})}{\usecounter{tmpabcd}}
\item \label{dist_alpha_cas_a} $\delta_0 \geq 2\cdot n!+1$,
\item \label{dist_alpha_cas_b} $\delta_1 \geq 4n$.
\end{list}
By~(\ref{ie_Qab}) we have
\begin{equation}
\ordz Q_{(a,b)}(\ullt{f}) \geq
\lfloor\frac{1}{n!}(a+1)(b+1)^n\rfloor > \frac{1}{n!}(a+1)(b+1)^n -1.
\end{equation}
In the case~\ref{dist_alpha_cas_a} one deduces
\begin{equation}
\ordz Q_{(a,b)}(\ullt{f}) \geq \frac{1}{n!}\left(\delta_0(\delta_1+1)^n-n!\right)
\end{equation}
and in order to show~(\ref{dist_alpha_est_ab}) it is sufficient to verify
\begin{equation} \label{dist_alpha_suffit_cas_a}
2\delta_0(\delta_1+1)^n-2\cdot n! \geq (\delta_0+1)(\delta_1+1)^n.
\end{equation}
The latter inequality is obvious for $\delta_0 \geq 2\cdot n!+1$ (and $\delta_1 \geq 0$).

In the case~\ref{dist_alpha_cas_b} the same procedure brings us to the point where it is sufficient to verify
(instead of~(\ref{dist_alpha_suffit_cas_a}))
\begin{equation} \label{dist_alpha_suffit_cas_b}
2(\delta_0+1)\delta_1^n-2\cdot n! \geq (\delta_0+1)(\delta_1+1)^n.
\end{equation}
We can rewrite this inequality as
\begin{equation} \label{ie_dist_alpha_suffit_cas_b}
\left(2\left(\frac{\delta_1}{\delta_1+1}\right)^n-1\right)(\delta_0+1) \geq \frac{2\cdot n!}{(\delta_1+1)^n}.
\end{equation}
 The l.h.s. of~(\ref{ie_dist_alpha_suffit_cas_b}) is an increasing function of $\delta_0$ and $\delta_1$, and the r.h.s. of~(\ref{ie_dist_alpha_suffit_cas_b}) is a decreasing function of $\delta_1$. So it is sufficient to verify this inequality for $\delta_0=0$ and $\delta_1=4n$. We can directly calculate
\begin{equation}
2\left(\frac{4n}{4n+1}\right)^n(0+1) > 1/2 > \frac{2\cdot n!}{(4n+1)^n},
\end{equation}
hence~(\ref{dist_alpha_suffit_cas_b}) is true for all the values $\delta_0 \geq 0$, $\delta_1 \geq 2n$. This completes the proof of~(\ref{dist_alpha_est_ab}).

We define
\begin{equation*}
\b{Q}(X_0,...,X_n)=Q_{(a,b)}(1,\b{z},X_0,...,X_n)^q,
\end{equation*}
where $q = \lceil 2\cdot n! \tilde{C} \rceil$; therefore we have
$\ordz Q_{(a,b)}(1,\b{z},1,\b{f}_1,...,\b{f}_n)^q \geq \tilde{C} (\delta_0+1)(\delta_1+1)^n$. As the coordinates of
$\ull{f}$ are algebraically independent over $\kk(\b{z})$, we have $\b{Q}(\ull{f}) \ne
0$.

It is easy to verify
\begin{equation*}
\begin{aligned}
&h(\b{Q})\leq\deg_{\ul{X}'}Q_{(a,b)} = a,\\
&\deg_{\ul{X}}\b{Q}=\deg_{\ul{X}}Q_{(a,b)} = b.
\end{aligned}
\end{equation*}
In the same time one obviously has
\begin{equation*}
\deg\b{Z} \geq 1
\end{equation*}
and, as $\b{Z}$ is  not defined over $\kk$ (see Remark~\ref{rem_LdT_Z_nonisotrivial}), one has
\begin{equation*}
h(\b{Z}) \geq 1.
\end{equation*}
In view of~(\ref{dist_alpha_choix_ab}) we obtain that the polynomial $\b{Q}$ makes the quantity~(\ref{dist_alpha_q}) strictly smaller than the minimum realized by $\b{P}_0$. So it can not vanish
on $\b{Z}$ (by the definition of $\b{P}_0$); in other words: $\b{Q}$ does not belong to $\I(\b{Z})$.

We apply Theorem~4.11 of chapter~3 of~\cite{NP} to
polynomial $\b{Q}(X_0,...,X_n)$
and to the ideal $\I(\b{Z})$ 0-dimensional over $\kk(\b{z})$.

Let $\ull{\alpha}\in\b{Z}$ realize the maximum of distance
from points of $\b{Z}$ to $\ull{f}$; in other words: let
$\Ord(\ull{f},\ull{\alpha})=\max_{\beta \in
\b{Z}}\Ord(\ull{f},\ull{\beta})$.

We define
\begin{equation} \label{def_theta}
\theta=\left\{ \begin{aligned} \Ordz \b{Q}(\ull{f}) &\mbox{, if }
\Ord(\ull{f},\ull{\alpha})>\Ordz \b{Q}(\ull{f})\\
\Ord_{\ull{f}} \I(\b{Z}) &\mbox{, if }
\Ord(\ull{f},\ull{\alpha}) \leq \Ordz \b{Q}(\ull{f})
\end{aligned}
\right.
\end{equation}

By Theorem~4.11 of chapter~3 of~\cite{NP} one has
\begin{equation} \label{thetaMaj}
\theta \leq h(\b{Q})\deg(\I(\b{Z}))+h(\I(\b{Z}))\deg(\b{Q})
\end{equation}
(in our case the base field is $\kk(\b{z})$ and all its valuations
are non-archimedean ones, so $\nu=0$ and the term $\nu m^2
\deg(\I(\b{Z})) \deg(\b{Q})$ is equal to zero in the
theorem in the reference).

We shall prove that the inequality
\begin{equation} \label{dist_alpha_ie1}
\Ord(\ull{f},\ull{\alpha}) \leq \Ordz \b{Q}(\ull{f})
\end{equation}
is in fact impossible.

Indeed, in this case $\theta=\Ord_{\ull{f}} \I(\b{Z})$, so~(\ref{thetaMaj}) implies
\begin{equation*} 
\Ord_{\ull{f}} \I(\b{Z}) \leq q \delta_0 \deg(\b{Z}) + q \delta_1 h(\b{Z}),
\end{equation*}
and we can weaken this inequality
\begin{equation*} 
\Ord_{\ull{f}} \I(\b{Z}) \leq \frac{q}{\min(\nu_0,\mu)}\left(\nu_0\delta_0 \deg(\b{Z}) + \nu_1\delta_1\deg(\b{Z}) + \mu\delta_1 h(\b{Z})\right).
\end{equation*}
Using the definition of $\delta_0$ and $\delta_1$ we deduce
\begin{equation} \label{majorationZf}
\begin{aligned}
\Ord_{\ull{f}} \I(\b{Z}) &\leq \frac{q}{\min(\nu_0,\mu)}\Big(\nu_0 h(\b{P}_0) \deg(\b{Z})\\&\qquad\qquad\qquad\qquad + \nu_1\deg\b{P}_0\deg(\b{Z}) + \mu\deg\b{P}_0h(\b{Z})\Big)\\
&\leq \frac{q}{\min(\nu_0,\mu)}\Big(\nu_0 h(\b{P}) \deg(\b{Z})\\&\qquad\qquad\qquad\qquad + \nu_1\deg\b{P}\deg(\b{Z}) + \mu\deg\b{P}h(\b{Z})\Big).
\end{aligned}
\end{equation}
Indeed, as $\b{P}$ vanishes on $\b{Z}$ one has
\begin{multline*}
 \nu_0 h(\b{P}_0) \deg(\b{Z}) + \nu_1\deg\b{P}_0\deg(\b{Z}) + \mu\deg\b{P}_0h(\b{Z}) \\ \leq \nu_0 h(\b{P}) \deg(\b{Z}) + \nu_1\deg\b{P}\deg(\b{Z}) + \mu\deg\b{P}h(\b{Z})
\end{multline*}
by the minimality defining $\b{P}_0$.
Then, applying~(\ref{LdTordZ}) (recall our notation $\b{Z}=\b{Z}_C(\b{P})$), one has
\begin{equation*}
\Ord_{\ull{f}}\I(\b{Z})> C^{\frac{1}{n}}c_n^{-\frac{n-1}{n}} \left(\nu_0 h(\b{P}) \deg(\b{Z}) + \nu_1\deg\b{P}\deg(\b{Z}) + \mu\deg\b{P}h(\b{Z})\right)
\end{equation*}
and gluing this inequality with~(\ref{majorationZf}) we obtain
\begin{multline*}
 C^{\frac{1}{n}}c_n^{-\frac{n-1}{n}} \left(\nu_0 h(\b{P}) \deg(\b{Z}) + \nu_1\deg\b{P}\deg(\b{Z}) + \mu\deg\b{P}h(\b{Z})\right)\\ < \frac{q}{\min(\nu_0,\mu)}\left(\nu_0 h(\b{P}) \deg(\b{Z}) + \nu_1\deg\b{P}\deg(\b{Z}) + \mu\deg\b{P}h(\b{Z})\right).
\end{multline*}
Simplifying $\nu_0 h(\b{P}) \deg(\b{Z}) + \nu_1\deg\b{P}\deg(\b{Z}) + \mu\deg\b{P}h(\b{Z})$ we deduce inequality
\begin{equation*}
 3\cdot n!\tilde{C}=C^{\frac{1}{n}}\min(\nu_0,\mu)c_n^{-\frac{n-1}{n}} < q = \lceil 2\cdot n!\tilde{C} \rceil
\end{equation*}
which contradicts the definition of~$q$ and $\tilde{C} \geq 1$ (recall that $\tilde{C}$ is defined at the end of the statement of this theorem and $\tilde{C} \geq 1$ in view of~(\ref{dist_alpha_Cestgrande})). So the inequality~(\ref{dist_alpha_ie1}) is impossible.

Unavoidably one has
\begin{equation*}
\Ord(\ull{f},\ull{\alpha})>\ordz \b{Q}(\ull{f}).
\end{equation*}
By construction of $\b{Q}$ one has $\ordz \b{Q}(\ull{f}) > \tilde{C} (\delta_0+1)(\delta_1+1)^n$, so we
deduce
\begin{equation}
\Ord(\ull{f},\ull{\alpha}) > \tilde{C} (\delta_0+1)(\delta_1+1)^n.
\end{equation}
\end{proof}

\section{Principal result}

In this section we introduce the main result of this paper in its full generality (a simplified version is given in Theorem~\ref{LMGP_simplified}) and prove it.

Recall that general framework imposed for this article is set up in subsection~\ref{subsection_general_framework}. So, we have an algebraically closed field $\kk$, a polynomial ring $\A=\kk[X_0',X_1',X_0,\dots,X_n]$, bi-graduated with respect to $\left(\deg_{\ul{X}'},\deg_{\ul{X}}\right)$, and a (set-theoretical) transformation $\phi:\A\rightarrow\A$ satisfying properties~(\ref{degphiQleqdegQ}) and~(\ref{condition_T2_facile}). In the statement below as well as forthcoming proofs we use various notions defined in subsections~\ref{definitions_comm_algebra} and~\ref{definitions_multiprojective_dg}. In particular, $m(I)$ (as well as $V_i$ and $e_{\phi}$) is defined in Definition~\ref{definDePP}, $i_0$ in Definition~\ref{def_i0} and $\ord_{\ul{f}}$ in Definition~\ref{defin_ord_xy}.
\begin{theorem}[Formal multiplicity lemma]\label{LMGP} Let $C_0 \in \mrr^+$ be a constant
such that for all bi-homogeneous prime ideal $\idq \subset \AnneauDePolynomes$ of rank $n$ one has
\begin{equation} \label{theoLMGP_condition_de_correctitude}
\ord_{\ull{f}}\idq \geq C_0 \Rightarrow \mbox{ the transformation $\phi$ is correct with respect to }\idq.
\end{equation}
Let $n_1\in\{1,\dots,n\}$ and $C_1\in\mrr^+$ and suppose that there exists a
constant $K_0 \in \mrr^{+}$ (depending only on $\phi$ and $\ullt{f}$) with the following property: for all equidimensional bi-homogeneous $\phi$-stable ideal $I\subset\AnneauDePolynomes$
 of rank $\geq n_1$ and such that $m(I)\leq (n+1)! \rho_{n+1}^{n+1}$,
and moreover all its associated prime ideals satisfy
\begin{equation} \label{theoLMGP_condition_ordp_geq_C0}
\ord_{\ull{f}}\idq \geq C_0,
\end{equation}
there exists a prime factor $\idq\in\Ass(\AnneauDePolynomes/I)$ that satisfies
\begin{equation} \label{RelMinN2}
\ord_{\ull{f}}(\idq) < K_0\left(\dd_{(0, n-\rg\idq+1)}(\idq)+\dd_{(1, n-\rg\idq)}(\idq)\right).
\end{equation}

Under these conditions, there exists a constant $K>0$ such that for all $P \in
\AnneauDePolynomes \setminus \{0\}$ satisfying for all $C\geq C_1$
\begin{equation}\label{condition_n1}
    i_0(\Z_C(P))\geq n_1
\end{equation}
(recall that the cycle $\Z_C(P)$ is introduced in Remark~\ref{genZP})
satisfy also
\begin{multline} \label{LdMpolynome2}
\ordz(P(\ullt{f})) \leq K\left((\mu+\nu_0)(\deg_{\ul{X}'}P+1)+\nu_1\deg_{\ul{X}}P\right)\\ \times\mu^{n-1}(\deg_{\ul{X}} P + 1)^n.
\end{multline}
\end{theorem}

\begin{remark} \label{rem_n1}
Condition~(\ref{condition_n1}) is always satisfied with $n_1=1$ and $C_1=0$ (in view of the definition of $i_0(Z(P))$, see Definition~\ref{def_i0} and Remark~\ref{rem_i0}).
Using this choice of $n_1$ and $C_0$ and forgetting also some restrictions on the $\phi$-stable ideals for which one has to provide~(\ref{theoLMGP_condition_ordp_geq_C0}) we obtain the statement of Theorem~\ref{LMGP_simplified}.
\end{remark}
\begin{remark}
Parameters $n_1$ and $C_1$ are introduced because in certain situations it is possible to provide direct lower estimate of $i_0(Z(P))$ better than 1 (see Remark~\ref{rem_i0}), so excluding the necessity of analysis of $\phi$-stable ideals of a small codimension. It could appear a decisive step, e.g. see proof of Proposition~4.11 in~\cite{EZ}.
\end{remark}



We shall deduce Theorem~\ref{LMGP} at the end of this section as a consequence of Proposition~\ref{PropositionLdMprincipal} and Lemma~\ref{LemmeProp13} below.

\begin{proposition} \label{PropositionLdMprincipal}
Let $P \in \AnneauDePolynomes$
and $C\in\mrr$ satisfy:
\begin{align}\label{C_bornee_enP}
C &< \frac{\ordz (P \circ \ullt{f}) - (\deg_{\ul{X}}P) \ordz (\ull{f})-\deg_{\ul{X}'}P}{n\left((\nu_0+\mu)(\deg_{\ul{X}'}P+1)+\nu_1\deg_{\ul{X}}P\right)\mu^{n-1}(\deg_{\ul{X}}P+1)^{n}}\\
C &\geq \left(\min(\nu_0,\mu)\right)^{-n}. \label{C_minore_numu}
\end{align}
Let $\idp$ be the ideal defined as $\idp=\I(\Z_C(P))$, where $\Z_C(P)$ is the cycle introduced in Remark~\ref{genZP}.
Suppose that for $i=i_0(\Z_C(P))$ one has
\begin{equation} \label{majoration_e_par_m}
e_{\phi}(V_i(\idp),\idp) \leq m(I_0(V_i(\idp),\idp));
\end{equation}
then
\begin{equation}\label{Cestsmall}
    C \leq \frac{(2\rho_{n+1})^nc_n^{n-1}}{\min(1;\lambda)^n\min(1;\mu)^n}.
\end{equation}
Moreover, for all polynomials $P\in\AnneauDePolynomes$, one has
\begin{equation} \label{estimationP}
\begin{aligned}
 \ordz &P(\ullt{f}(\b{z}))\leq\max\left(\frac{n}{\left(\min(\nu_0,\mu)\right)^{n}},\frac{(2\rho_{n+1})^nc_n^{n-1}}{\min(1;\lambda)^n\min(1;\mu)^n}\right)\\
 &\times\left((\mu+\nu_0)(\deg_{\ul{X}'}P+1)+\nu_1\deg_{\ul{X}}P\right)\mu^{n-1}(\deg_{\ul{X}} P+1)^n\\
 &+(\ordz \ull{f})(\deg_{\ul{X}} P)+\deg_{\ul{X}'}P.
\end{aligned}
\end{equation}
\end{proposition}
\begin{proof}
Note that for $\deg_{\ul{X}}P=0$ the conclusion of the proposition 
is automatically satisfied. Thus we need only to treat the case $\deg_{\ul{X}}P\geq 1$.

{\it Ad absurdum} assume
\begin{equation} \label{cnstCgrande}
 C > \frac{(2\rho_{n+1})^nc_n^{n-1}}{\min(1;\lambda)^n\min(1;\mu)^n}.
\end{equation}

Recall that $i_0(\Z_C(P))\geq 1$ is
the largest index $i \in \{1,...,n\}$ such that $\rg(V_i\A_{\idp}) \geq i$ (see Definition~\ref{def_i0}).
We put $e_0$ the largest integer $\leq e_{\phi}(V_{i_0},\idp)$
such that $V_{i_0}+...+{\phi}^{e_0}(V_{i_0})\subset\idp$ (we use the notation $V_{i_0}$ as a shorthand for $V_{i_0}(\idp)$). Note that the assumption~(\ref{majoration_e_par_m}) implies that
$e_{\phi}(V_{i_0},\idp)$ is finite, so $e_0$ is a well-defined integer.

\smallskip

Let $Q$ be a generator of ${\phi}^{e_0}(V_{i_0})$; by Lemma~\ref{majorationphinQ} one has
\begin{equation} \label{PropositionLdMprincipal_majoration_deg_Q}
\begin{aligned}
&\deg_{\ul{X}}Q \leq \mu^{e_0}\rho_{i_0}\delta_1(\idp),\\
&\deg_{\ul{X}'}Q \leq (\nu_0\delta_0(\idp)+e_0\nu_1\delta_1(\idp))\max(\nu_0,\mu)^{e_0-1}\rho_{i_0}.
\end{aligned}
\end{equation}

With the substitution $(X_0':X_1')=(1:\b{z})$ we can consider $Q$ as a polynomial $\b{Q}$ of $\kk[\b{z}][X_0:...:X_n]$. We define $\b{Z}=\b{Z}_C(P)$.

Let $\b{\alpha}\in\b{Z}$. By Lemma~\ref{Representants}, b), there is a system of projective coordinates $\ull{\alpha}$ satisfying
\begin{eqnarray*}
 \ordz\ull{\alpha} &=& \ordz\ull{f}, \\
  \ordz(\ull{\alpha} - \ull{f}) - \ordz(\ull{f}) &=& \Ordz(\b{\alpha},\b{f}).
\end{eqnarray*}
In view of $\ordz\ull{f}=0$, we deduce immediately
\begin{eqnarray}
  \ordz\ull{\alpha} &=& 0,\label{PropositionLdMprincipal_alpha_c1} \\
  \ordz(\ull{\alpha} - \ull{f}) &=& \ordz(\ull{\alpha}\wedge\ull{f}). \label{PropositionLdMprincipal_alpha_c2}
\end{eqnarray}

We fix a choice of projective coordinate systems satisfying~(\ref{PropositionLdMprincipal_alpha_c1}) and~(\ref{PropositionLdMprincipal_alpha_c2})
for all $\b{\alpha}\in\b{Z}$.

We claim that
\begin{equation} \label{point_Crucial}
 \begin{aligned}
  \ordz(\b{\phi}(\b{Q})(\ull{\alpha}))
	&\geq \min(\ordz(\b{\phi}(\b{Q})(\ull{f})),\ordz(\ull{\alpha}\wedge\ull{f})).
 \end{aligned}
\end{equation}
Indeed,
\begin{equation*}
\begin{aligned}
    \ordz&\left(\b{\phi}(\b{Q})(\ull{\alpha})\right)=\ordz\left(\left(\b{\phi}(\b{Q})(\ull{\alpha})-\b{\phi}(\b{Q})(\ull{f})\right)+\b{\phi}(\b{Q})(\ull{f})\right)\\
    &\geq\min\left(\ordz\left(\b{\phi}(\b{Q})(\ull{\alpha})-\b{\phi}(\b{Q})(\ull{f})\right),\ordz\left(\b{\phi}(\b{Q})(\ull{f})\right)\right)\\
    &\geq\min\left(\ordz\left(\ull{\alpha}-\ull{f}\right),\ordz\left(\b{\phi}(\b{Q})(\ull{f})\right)\right)\\
    &\geq\min\left(\ordz\left(\ull{\alpha}\wedge\ull{f}\right),\ordz\left(\b{\phi}(\b{Q})(\ull{f})\right)\right).
\end{aligned}
\end{equation*}
Then, by~(\ref{condition_T2_facile}) and in view of $\b{Q}(\ull{\alpha})=0$ (according to the choice of $e_0$),
\begin{equation} \label{point_Crucial_1}
 \begin{split}
  \ordz\left(\b{\phi}(\b{Q})(\ull{f})\right) &\geq \lambda \ordz\b{Q}(\ull{f})\\
	&\geq \lambda \ordz\big(\b{Q}(\ull{f})-\b{Q}(\ull{\alpha})\big)\\
	&\geq \lambda\ordz\left(\ull{\alpha}\wedge\ull{f}\right).
 \end{split}
\end{equation}
We deduce from~(\ref{point_Crucial}) and~(\ref{point_Crucial_1})
\begin{equation} \label{PropositionLdMprincipal_ie0}
 \ordz(\b{\phi}(\b{Q})(\ull{\alpha})) \geq \min(1,\lambda) \ordz(\ull{\alpha}\wedge\ull{f})
\end{equation}
for all $\b{\alpha}\in\b{Z}$.

By~(\ref{PropositionLdMprincipal_ie0}) one has
\begin{equation} \label{PropositionLdMprincipal_ie1_1}
 \begin{split}
  &\frac{1}{\min(1,\lambda)}\sum_{\ull{\alpha} \in \b{Z}_{C}(P)}\left(\ordz\left(\phi(\b{Q})(\ull{\alpha})\right)\right)
  \\&\geq\sum_{\ull{\alpha} \in \b{Z}_{C}(P)}\left(\ordz(\ull{\alpha}\wedge\ull{f})\right)=:M
 \end{split}
\end{equation}
(note that $M$ is equal the l.h.s. of~(\ref{LdTordZ})). By definition of $\b{Z}_C(P)$ (see Definition~\ref{defZP} and Remark~\ref{genZP}) and with~(\ref{cnstCgrande}) we estimate
\begin{multline} \label{PropositionLdMprincipal_ie1}
  M>C^{\frac{1}{n}}c_n^{-\frac{n-1}{n}}\left(\nu_0\deg(\b{Z})\deg_{\b{z}}P+\nu_1\deg(\b{Z})\deg_{\ul{X}}P + \mu h(\b{Z})\deg_{\ul{X}}P\right)\\
	\geq \frac{2\rho_{n+1}}{\min(1,\lambda)\min(1,\mu)}\Big(\nu_0\deg(\b{Z})\deg_{\b{z}}P\\ + \nu_1\deg(\b{Z})\deg_{\ul{X}}P + \mu h(\b{Z})\deg_{\ul{X}}P\Big).
\end{multline}
We deduce from~(\ref{PropositionLdMprincipal_ie1_1}), (\ref{PropositionLdMprincipal_ie1}) and taking into account~(\ref{cnstCgrande})
\begin{multline} \label{PropositionLdMprincipal_ie15}
  \sum_{\ull{\beta} \in \b{Z}_{C}(P)}\ordz\left(\phi(\b{Q})(\ul{\beta})\right)
	>\frac{2\rho_{n+1}}{\min(1,\mu)}\\ \times\left(\nu_0\deg(\b{Z})\deg_{\b{z}}P + \nu_1\deg(\b{Z})\deg_{\ul{X}}P + \mu h(\b{Z})\deg_{\ul{X}}P\right).
\end{multline}
Also Liouville's inequality~(\ref{iet_main})
implies (if $\b{\phi}(\b{Q})$ does not vanish on $\b{Z}_{C}(P)$)
\begin{equation} \label{PropositionLdMprincipal_ie2}
 \begin{split}
  &\sum_{\ull{\beta} \in \b{Z}_{C}(P)}\ordz\left(\phi(\b{Q})(\ul{\beta})\right)\leq\deg(\b{Z})h(\phi(\b{Q})) + h(\b{Z})\deg\phi(\b{Q})\\
  &\leq \max(\mu,\nu_0)^{e_0}\rho_{i_0}\left(\nu_0\deg(\b{Z})\delta_0 + \nu_1(e_0+1)\deg(\b{Z})\delta_1 + \mu h(\b{Z})\delta_1\right)\\
  &\leq \max(\mu,\nu_0)^{e_0}(e_0+1)\rho_{i_0}\left(\nu_0\deg(\b{Z})\delta_0 + \nu_1\deg(\b{Z})\delta_1 + \mu h(\b{Z})\delta_1\right)
 \end{split}
\end{equation}
(the second inequality here comes from~(\ref{PropositionLdMprincipal_majoration_deg_Q})).

According to the definition of $e_0$, the hypothesis~(\ref{majoration_e_par_m}) and Lemma~\ref{LemmeCor14NumberW} we have
\begin{equation} \label{PropLdMprincipal_ie_e0_leq_m}
e_0 \leq e_{\phi}(V_{i_0},\idp) \leq m(I_0(V_{i_0},\idp)) \leq \nu(n+1)!\rho_{i_0}^{n+1},
\end{equation}
and so $\max(\mu,\nu_0)^{e_0}(e_0+1)\rho_{i_0}\leq \rho_{i_0+1} \leq \rho_{n+1}$ by definition of $\rho_{n+1}$.
Thus~(\ref{PropositionLdMprincipal_ie15}) and~(\ref{PropositionLdMprincipal_ie2}) lead to:
\begin{multline*}
    \frac{2\rho_{n+1}}{\min(1,\mu)}\left(\nu_0\deg(\b{Z})\deg_{\b{z}}P + \nu_1\deg(\b{Z})\deg_{\ul{X}}P + \mu h(\b{Z})\deg_{\ul{X}}P\right)\\
    <\rho_{n+1}\left(\nu_0\deg(\b{Z})\delta_0 + \nu_1\deg(\b{Z})\delta_1 + \mu h(\b{Z})\delta_1\right).
\end{multline*}
This inequality contradicts Definition~\ref{def_delta}, thus $\b{\phi}(\b{Q})(\ull{\alpha})=0$.

So we have
\begin{equation} \label{PropositionLdMprincipal_incl_e0p1}
{\phi}^{e_0+1}(V_{i_0}) \subset \idp,
\end{equation}
and this inclusion contradicts the definition of $e_0$ if $e_0 < e_{\phi}(V_{i_0},\idp)$. We conclude
$e_0=e_{\phi}(V_{i_0},\idp)$.

Moreover, (\ref{PropositionLdMprincipal_incl_e0p1}) implies
\begin{equation} \label{PropositionLdMprincipal_est_rg_rgP}
\rg\left((V_{i_0}+...+\phi^{e_0+1}(V_{i_0}))\AnneauDePolynomes_{\idp}\right)\leq\rg\left(\idp\AnneauDePolynomes_{\idp}\right)=n.
\end{equation}
As $e_0+1>e_{\phi}(V_{i_0},\idp)$ and by definition of $e_{\phi}(V_{i_0},\idp)$ we have
\begin{equation} \label{e0plus1Total}
\begin{split}
 \rg\left((V_{i_0}+...+\phi^{e_0+1}(V_{i_0}))\AnneauDePolynomes_{\idp}\right) > \rg(V_{i_0}\AnneauDePolynomes_\idp) \geq i_0,
\end{split}
\end{equation}
we obtain
\begin{equation*}
 \rg(V_{i_0+1}\AnneauDePolynomes_\idp)\geq\rg\left((V_{i_0}+...+\phi^{e_0+1}(V_{i_0}))\AnneauDePolynomes_{\idp}\right)\geq i_0+1.
\end{equation*}
If $i_0<n$ this inequality contradicts the definition of $i_0$, and if $i_0=n$ inequality~(\ref{e0plus1Total})
implies
\begin{equation*}
\rg\left((V_{i_0}+...+\phi^{e_0+1}(V_{i_0}))\AnneauDePolynomes_{\idp}\right) > n,
\end{equation*}
in contradiction with~(\ref{PropositionLdMprincipal_est_rg_rgP}).

So, we have verified that the hypothesis~(\ref{cnstCgrande}) can not being satisfied, establishing therefore~(\ref{Cestsmall}).

It remains to verify~(\ref{estimationP}). We fix an arbitrary polynomial $P\in\A$ and consider the set $\mathcal{M}(P)$ of reals $C$ satisfying~(\ref{C_bornee_enP}) and~(\ref{C_minore_numu}) (for our chosen and fixed polynomial $P$).
If $\mathcal{M}(P)=\emptyset$, we have
\begin{equation*}
    \left(\min(\nu_0,\mu)\right)^{-n}\geq\frac{\ordz (P \circ \ullt{f}) - (\deg_{\ul{X}}P) \ordz (\ull{f})-\deg_{\ul{X}'}P}{n\left((\mu+\nu_0)(\deg_{\ul{X}'}P+1)+\nu_1\deg_{\ul{X}}P\right)\mu^{n-1}(\deg_{\ul{X}}P+1)^{n}}
\end{equation*}
obtaining immediately~(\ref{estimationP}).

In the opposite case, if $\mathcal{M}(P)\ne\emptyset$, we let $C_s$ denote the upper bound of $\mathcal{M}(P)$, which is a real finite number: the inequality~(\ref{C_bornee_enP})
shows 
\begin{equation*}
    C_s=\frac{\ordz (P \circ \ullt{f}) - (\deg_{\ul{X}}P) \ordz (\ull{f})-\deg_{\ul{X}'}P}{n\left((\mu+\nu_0)(\deg_{\ul{X}'}P+1)+\nu_1\deg_{\ul{X}}P\right)\mu^{n-1}(\deg_{\ul{X}}P+1)^{n}}.
\end{equation*}
In the first part of the proof we have established the inequality~(\ref{Cestsmall}) for all the elements of $\mathcal{M}(P)$, therefore $C_s$ also satisfies this
inequality, hence~(\ref{estimationP}).
\end{proof}

\begin{lemma}\label{LemmeProp13}
Let $\idp$ be a prime ideal of $\A$
such that the transformation $\phi$ is correct with respect to this ideal and
let $V \subset \A$, $V\ne\{0\}$, be a $\kk$-linear subspace of
$\A$. If $e_{\phi}(V,\idp)>m(I_0(V,\idp))$, then there
exists an equidimensional $\phi$-stable ideal $J$ such that
\begin{list}{\alph{tmpabcd})}{\usecounter{tmpabcd}}
 \item \label{LemmeProp13_a} $I_0(V,\idp) \subset J \subset \idp$, 
 \item \label{LemmeProp13_b} $\rg(J)=\rg(I_0(V,\idp))$, 
 \item \label{LemmeProp13_c} all the primes associated to $J$ are contained in $\idp$. 
\end{list}
In particular,
\begin{equation} \label{LemmeProp13_en_part}
\begin{aligned}
    &m(J) \leq m(I_0(V,\idp)),\\
    &\dd_{(1, n-\rg J)} J \leq \dd_{(1, n-\rg I_0(V,\idp))}(I_0(V,\idp)),\\
    &\dd_{(0, n-\rg J+1)} J \leq \dd_{(0, n-\rg I_0(V,\idp)+1)}(I_0(V,\idp)).
\end{aligned}
\end{equation}
\end{lemma}
\begin{proof}
We shall use the shorthand $m$ for $m(I_0(V,\idp))$, and for each $e=0,1,2,...$ we use the notation $\tilde{J_e}:=I_e(V,\idp)$, $J_e:=\eq(\tilde{J}_e)$. In view of these definitions we have
\begin{equation}
\tilde{J}_e \subset \tilde{J}_{e+1} \subset J_{e+1} \qquad e=0,1,2,....
\end{equation}
and
\begin{equation} \label{LemmeProp13_rkJep1_rkJe}
    \rk J_e=\rk\tilde{J}_e\leq\rk\tilde{J}_{e+1}=\rk J_{e+1}.
\end{equation}

Note that if $\rg(\tilde{J}_{e+1})=\rg(\tilde{J}_e)$, that is to say if $\rg(J_{e+1})=\rg(J_e)$, then
\begin{equation}\label{Jep1_sps_Je}
    J_{e+1} \supset J_e
\end{equation}

(because obviously we always have the inclusion $\tilde{J}_{e+1} \supset \tilde{J_e}$).

If
\begin{equation}\label{phiep1_nss_Je}
    \phi^{e+1}(V) \nsubseteq J_e,
\end{equation}
then either:
\begin{equation}\label{ie_Jep1_Je}
    \rg(J_{e+1}) > \rg(J_e)
\end{equation}
or
\begin{equation}\label{ie_mJep1_mJe}
    \rg(J_{e+1})=\rg(J_e) \text{ and } m(J_{e+1}) < m(J_e).
\end{equation}
Indeed, (\ref{phiep1_nss_Je})) implies that there is an element $x\in J_{e+1}\setminus J_e$
(because $J_{e+1} \supset \tilde{J}_{e+1} \supset \phi^{e+1}(V)$).
\emph{A fortiori} either this $x$ is not a zero divisor in $\AnneauDePolynomes/J_e$ (and as a consequence $\rg(J_{e+1})>\rg(J_e)$, i.e. we have~(\ref{ie_Jep1_Je})),
or, if we assume that $x$ \emph{is} a zero divisor in $\AnneauDePolynomes/J_e$, then $x$ does not belong to at least one primary component of the equidimensional ideal $J_e$.
In the latter case we can also suppose $\rg(J_{e+1}) = \rg(J_e)$ (otherwise we have~(\ref{ie_Jep1_Je}) and the claim is established), so we also have $J_e\subset J_{e+1}$ (see~(\ref{Jep1_sps_Je}))
and each primary component of $J_{e+1}$ contains at least one primary component of $J_e$
(because otherwise by localization we couldn't obtain $J_e\subset J_{e+1}$).
So for all $\idq\subset\Spec(\AnneauDePolynomes)$ we have
$l_{\AnneauDePolynomes_{\idq}}((\AnneauDePolynomes/J_{e+1})_{\idq})\leq l_{\AnneauDePolynomes_{\idq}}((\AnneauDePolynomes/J_{e})_{\idq})$
and in view of~(\ref{definDePP_defin_m}) we conclude $m(J_{e+1})\leq m(J_e)$.
But there exists an element $x\in J_{e+1}$ that does not belong to at least one primary component $\idq$ of $J_e$, we denote the radical of this primary component $\idq$ by $\idp$.
So we obtain $l_{\AnneauDePolynomes_{\idp}}((\AnneauDePolynomes/J_{e+1})_{\idp})<l_{\AnneauDePolynomes_{\idp}}((\AnneauDePolynomes/J_e)_{\idp})$, and
thus $m(J_{e+1})<m(J_e)$.

By the definition of $e_{\phi}(V,\idp)$ (see Definition~\ref{definDePP}) we have
\begin{equation} \label{rgrgJJ}
 \rg\left(J_e\right)=\rg\left(J_0\right) \mbox{ pour tout } e=0,...,e_{\phi}(V,\idp).
\end{equation}
As $e_{\phi}(V,\idp)\geq m$, if $\phi^{e+1}(V) \nsubseteq J_e$ for $e=0,...,m-1$, we have
\begin{equation} \label{ordonneJJ}
 1 \leq m(J_m) < m(J_{m-1}) < \dots < m(J_0)=m,
\end{equation}
thus we have $m+1$ pairwise distinct integers $m(J_0), \dots , m(J_m)$ between $1$ and $m$, this situation been impossible
we conclude that there exists $0 \leq e_0 \leq m-1$ such that
\begin{equation}
 \phi^{e_0+1}(V) \subset J_{e_0}.
\end{equation}
Thus we have (in view of the definition of the ideal $\tilde{J}_{e_0}$)
\begin{equation} \label{Jep1_eq_Je}
\phi(\tilde{J}_{e_0})\subset \left(\tilde{J}_{e_0}+\phi^{e_0+1}(V)\right)\AnneauDePolynomes_{\idp}\cap\AnneauDePolynomes\subset J_{e_0}+J_{e_0}\subset J_{e_0},
\end{equation}
because $\tilde{J}_{e_0}\subset J_{e_0}$ by definition of $J_{e_0}$.

As $\phi$ is correct with respect to $\idp$ (and by construction of $J_{e_0}$ all its associated primes are contained in $\idp$),
(\ref{Jep1_eq_Je}) implies
\begin{equation*}
    \phi(J_{e_0})\subset J_{e_0},
\end{equation*}
so the equidimensional ideal $J_{e_0}$ is $\phi$-stable.

\smallskip

We are going to show that the choice $J=J_{e_0}$ satisfies all the conditions of the lemma.

The construction of $J_{e_0}$ shows that all its associated  primes
are contained in $\idp$, it gives us the point~\ref{LemmeProp13_c}. Further, we deduce the point~\ref{LemmeProp13_a} with the evident remark $I_0(V,\idp)=\tilde{J}_0 \subset J_{e_0}$.

The point \ref{LemmeProp13_b} is implied by~(\ref{rgrgJJ}).

Finally, (\ref{LemmeProp13_en_part}) is a consequence of a), b) and c). So the ideal $J=J_{e_0}$ satisfies all the properties from the statement of the proposition.
\end{proof}

\begin{proof}[Proof of Theorem~\ref{LMGP}.] We define
\begin{multline} \label{Cestgrande}
 C=1+\max\Bigg(c_n^{n-1}C_0^{n}, \left(\min(\nu_0,\mu)\right)^{-n}, \Ciso,\\ \left(n!c_n\left(1+\frac{\nu_1}{\max(\mu,\nu_0)}\right)\rho_n^{n}K_0\right)^{n},C_1\Bigg)
\end{multline}
(recall that $c_n$ is defined in~(\ref{def_cn}), $C_0$, $C_1$, $K_0$ come from the statement of Theorem~\ref{LMGP}
(assumptions~(\ref{theoLMGP_condition_ordp_geq_C0}), (\ref{condition_n1})),
$\Ciso$ is provided by the remark~\ref{rem_LdT_Z_nonisotrivial}). Let $P \in \AnneauDePolynomes\setminus\{0\}$ be a polynomial that does not satisfy~(\ref{LdMpolynome2}) for
\begin{equation}\label{preuve_theo_defK}
    K=\max\left(2nC,\left(\frac{2\rho_{n+1}c_n}{\max(1,\lambda)\max(1,\mu)}\right)^n\right).
\end{equation}
Then it satisfies~(\ref{ordPplusqueb2}). In particular, $C$ and $P$ satisfy
conditions~(\ref{C_bornee_enP}) and~(\ref{C_minore_numu}) of Proposition~(\ref{PropositionLdMprincipal}).

In the sequel we use the notation $\idp=\I(\Z_C(P))$, where $\Z_C(P)$ is the cycle introduced in the remark~\ref{genZP}.
In view of~(\ref{LdTordZ}) and~(\ref{Cestgrande}), we have $\ord_{\ullt{f}}\idp > C_0$ and thus, using the assumption~(\ref{theoLMGP_condition_de_correctitude}),
$\phi$ is correct with respect to $\idp$. Moreover, $\Z_C(P)$ is projected onto $\mpp^1$ (cf. the remark~\ref{rem_LdT_Z_nonisotrivial}).

Recall that $V_i=V_i(\Z_C(P))$ is introduced in the definition~(\ref{V_irho_i}) and
$e_{\phi}$, $m$ are introduced
in Definition~\ref{definDePP}: (\ref{definDePP_defin_e}) and (\ref{definDePP_defin_m}) respectively
If for $i=i_0(\Z_C(P))$ we have
\begin{equation}\label{e_bornee}
    e_{\phi}(V_i,\idp) \leq m(I_0(V_{i},\idp)),
\end{equation}
we verify~(\ref{majoration_e_par_m}) and we can apply Proposition~\ref{PropositionLdMprincipal},
it gives us~(\ref{LdMpolynome2}) in view of our choice of~$K$. This contradicts our hypothesis that $P$ does not satisfy~(\ref{LdMpolynome2}).
On the other hand, if~(\ref{e_bornee}) is not satisfied, we apply Lemma~\ref{LemmeProp13} to the ideal $\idp$
and the vectorial space $V=V_i(\idp)$ (we recall the notation $i=i_0(\Z_C(P))$).

We denote by $J$ the equidimensional $\phi$-stable ideal provided by
Lemma~\ref{LemmeProp13}. In view of the property~$b)$ of this proposition we have $\rg J=\rg\left(V_i\right)\geq i\geq n_1$.

We verify (In view of~(\ref{LemmeProp13_en_part}))
\begin{equation}
\begin{aligned}
& m(J) \leq m(I_0(V,\idp)),\\
& \dd_{(0,n-\rg J+1)}(J) \leq \dd_{(0,n-\rg J+1)}(I_0(V_i,\idp)).
\end{aligned}
\end{equation}
As $\V(\idp)=Z_C(P)$ is projected onto $\mpp^1$ we have by Lemma~\ref{LemmeCor14NumberW}
$$
m(I_0(V,\idp)) \leq \nu(n+1)!\rho_i^{n+1}
$$
and also $\delta_1(\idp) \geq 1$.

Recall that $I_0(V_i,\idp)\subset\idp$ and thus
\begin{equation*}
\rg I_0(V_i,\idp) \leq \rg\idp = n.
\end{equation*}
As the ideal~$I_0(V_i,\idp) \subset \idp$ is extended-contracted of an ideal generated by polynomials of bi-degree
$\leq \left(\rho_i\left(\delta_0(\idp)+\frac{\nu_1}{\max(\mu,\nu_0)}\delta_1(\idp)\right), \rho_i\delta_1(\idp)\right)$, we have
\begin{equation} \label{est_basique1}
\begin{aligned}
\dd_{(0, n-\rg I_0(V_i,\idp)+1)}(I_0(V_i,\idp)) &\leq \left(\delta_0(\idp)+\frac{\nu_1}{\max(\mu,\nu_0)}\delta_1(\idp)\right)\\&\qquad\times\delta_1(\idp)^{\rg I_0(V_i,\idp)-1}\rho_i^{\rg I_0(V_i,\idp)}\\
&\leq\left(1+\frac{\nu_1}{\max(\mu,\nu_0)}\right)\rho_n^{n}\\&\times\left(\delta_0(\idp)+1\right)\left(\delta_1(\idp)+1\right)^{\rg I_0(V_i,\idp)},
\end{aligned}
\end{equation}
and finally
\begin{equation} \label{estdegW}
 \dd_{(0, n-\rg J+1)}J \leq \left(1+\frac{\nu_1}{\max(\mu,\nu_0)}\right)\rho_n^{n}\left(\delta_0(\idp)+1\right)\left(\delta_1(\idp)+1\right)^{\rg J}.
\end{equation}
As $J$ is an equidimensional ideal, we obtain for all $\idq\in\Ass(\AnneauDePolynomes/J)$
\begin{equation} \label{estdegQ}
\begin{aligned}
 \dd_{(0, n-\rg\idq+1)}\idq&\leq\dd_{(0, n-\rg J+1)}J\\
 &\leq\left(1+\frac{\nu_1}{\max(\mu,\nu_0)}\right)\rho_n^{n}\left(\delta_0(\idp)+1\right)\left(\delta_1(\idp)+1\right)^{\rg\idq}.
\end{aligned}
\end{equation}
The same calculation for $\dd_{(1, n-\rg\idq)}\idq$ gives us
\begin{equation} \label{estdeg1Q}
\begin{aligned}
 \dd_{(1, n-\rg\idq)}\idq
 \leq\rho_n^{n}\left(\delta_1(\idp)+1\right)^{\rg\idq}.
\end{aligned}
\end{equation}

As $P$ and $C$ satisfy~(\ref{ordPplusqueb2}), by Lemma~\ref{dist_alpha} there exists a point $\ull{\alpha}\in\b{Z}_C(P)$
satisfying~(\ref{Cdirect}) with $\tilde{C} = \frac{C^{\frac{1}{n}}\min(\nu_0,\mu)}{3 n!c_n}\geq\left(1+\frac{\nu_1}{\max(\mu,\nu_0)}\right)\rho_n^{n}K_0$
(the last inequality is implied by the definition~(\ref{Cestgrande})),
and thus one has for all $\idq\in\Ass(\AnneauDePolynomes/J)$ (in view of Lemma~\ref{LemmeProp13}, point~\ref{LemmeProp13_c})
\begin{equation} \label{estordQ}
 \ord_{\ull{f}}\idq \geq \ord(\ull{f},\ull{\alpha}) > \left(1+\frac{\nu_1}{\max(\mu,\nu_0)}\right)\rho_n^{n} K_0(\delta_0(\idp)+1)(\delta_1(\idp)+1)^n.
\end{equation}

The estimates~(\ref{estdegQ}), (\ref{estdeg1Q})
and~(\ref{estordQ}) put together (and been verified for \emph{all} $\idq\in\Ass(\AnneauDePolynomes/J)$)
contradict~(\ref{RelMinN2}). So, the assumption~(\ref{ordPplusqueb2})
with $C$ given by~(\ref{Cestgrande}) is untenable and we deduce~(\ref{LdMpolynome2}) with our choice of~$K$.
It contradicts again our assumption that $P$ does not satisfy~(\ref{LdMpolynome2}).

Finally, we conclude that polynomial~$P$ does not satisfying~(\ref{LdMpolynome2}) can not exist, and this completes the proof.
\end{proof}

\section{Applications}

\subsection{Multiplicity estimates for solutions
of algebraic differential equations} \label{subsection_ApplicationsDifferential}

In this subsection we consider an $n$-tuple $\ull{f}=(f_1(\b{z}),\dots,f_n(\b{z}))$ of analytic functions (or, more generally, power series) satisfying the system of differential equations
\begin{equation} \label{syst_diff}
f_i'(\b{z})=\frac{A_i(\b{z},\ull{f})}{A_0(\b{z},\ull{f})}, \quad i=1,\dots,n,
\end{equation}
where $A_i(\b{z},X_1,\dots,X_n)\in\kk[\b{z},X_1,\dots,X_n]$ for $i=0,...,n$ (we suppose that $A_0$ is a non-zero polynomial).

We associate to the system~(\ref{syst_diff}) the following differential operator
\begin{equation} \label{defD}
D = A_0(\b{z}, X_1,\dots, X_n)\diff{\b{z}} + \sum_{i=1}^nA_i(\b{z}, X_1,\dots, X_n)\diff{X_i}.
\end{equation}
This operator is an application $D:\kk[\b{z},X_1,\dots,X_n]\rightarrow\kk[\b{z},X_1,\dots,X_n]$. We also consider $D$ as acting on $\AnneauDePolynomes=\kk[X_0',X_1'][X_1,\dots,X_n]$ defining
\begin{equation} \label{defhD}
D = {^h\!A}_0(X_0',X_1', X_1,\dots, X_n)\diff{X_1'}+\sum_{i=1}^n{^h
\!A_i}(X_0',X_1', X_1,\dots, X_n)\diff{X_i},
\end{equation}
where $^h\!P$ denotes the bi-homogenization of the polynomial $P\in
\kk[\b{z},X_1,\dots,X_n]$:
\begin{equation*}
  ^h\!P(X_0',X_1', X_1,\dots, X_n):=X_0'^{\deg_{\b{z}}P}\cdot
X_0^{\deg_{\ul{X}}P}\cdot P\left(\frac{X_1'}{X_0'},\frac{X_1}{X_0},
\dots,\frac{X_n}{X_0}\right).
\end{equation*}
One readily verifies~$D({^h\!P})={^h\!\left(D(P)
\right)}$, so the application $D:\AnneauDePolynomes\rightarrow
\AnneauDePolynomes$ is exactly the "bi-homogenization" of $D:
\kk[\b{z},X_1,\dots,X_n]\rightarrow\kk[\b{z},X_1,\dots,X_n]$.

The application $D$ is a correct application with respect to any
ideal $\idp \subset \AnneauDePolynomes$, according to Corollary~\ref{exemple_OpDiff_Correcte}.

We are going to deduce from Theorem~\ref{LMGP} an improvement of the following Nesterenko's famous theorem (proved in~\cite{N1996}):

\begin{theorem}[Nesterenko, see Theorem~1.1 of Chapter~10, \cite{NP}] \label{theoNesterenko_classique}
Suppose that functions
\begin{equation*}
\ull{f} = (f_1(\b{z}),\dots,f_n(\b{z})) \in \mcc[[\b{z}]]^n
\end{equation*}
are analytic at the point $\b{z}=0$ and form a solution of the system~(\ref{syst_diff}) with $\kk=\mcc$.
If there exists a constant $K_0$ such that every $D$-stable prime ideal $\idp \subset \mcc[X_1',X_1,\dots,X_n]$,
$\idp\ne(0)$, satisfies
\begin{equation} \label{ordIleqKdegI}
\min_{P \in \idp}\ordz P(\b{z},\ull{f}) \leq K_0,
\end{equation}
then there exists a constant $K_1>0$ such that for any polynomial $P \in
\mcc[X_1',X_1,\dots,X_n]$, $P\ne 0$, the following inequality holds
\begin{equation} 
\ordz(P(\b{z},\ull{f})) \leq K_1(\deg_{\ul{X}'} P + 1)(\deg_{\ul{X}} P + 1)^n.
\end{equation}
\end{theorem}

\begin{remark}
Assuming $A_0(0,\ull{f}(0))\ne 0$ in the system~(\ref{syst_diff}), it is easy to verify the condition~(\ref{ordIleqKdegI}), cf.~\cite{NP}, chapitre 10, example~1 (p.~150).
Also, the condition~(\ref{ordIleqKdegI}) is established in the case when the polynomials $A_i$, $i=0,\dots,n$ are of degree 1 in $X_1,\dots,X_n$, cf.~\cite{N1974}.
In the latter case the proof is based on the differential Galois theory.
\end{remark}

Using Theorem~\ref{LMGP} we can replace~(\ref{ordIleqKdegI}) in Theorem~\ref{theoNesterenko_classique} by
a weaker assumption~(see~(\ref{ordIleqKdegI2}) below and~(\ref{RelMinN2})). In the same time our result is valid with the base field of an arbitrary characteristic.

\begin{theorem} \label{LMGPD}
Let
\begin{equation*}
\ull{f} = (f_1(\b{z}),\dots,f_n(\b{z})) \in \kk[[\b{z}]]^n
\end{equation*}
be a set of formal power series forming a solution of the system~(\ref{syst_diff}) and let $n_1\in\{1,\dots,n\}$, $C_1\in\mrr^+$ be two numbers.
We suppose that there is a constant $K_0$ such that all bi-homogeneous primary $D$-stable ideals
$\idp \subset \AnneauDePolynomes$ of rank $\geq n_1$ satisfy
\begin{equation} \label{ordIleqKdegI2}
\ord_{\ull{f}}\idp \leq K_0\left(\dd_{(0,n-\rg\idp+1)}\idp+\dd_{(1,n-\rg\idp)}\idp\right),
\end{equation}
or, if $\car\kk=0$, we suppose that inequality~(\ref{ordIleqKdegI2}) holds only for bi-homogeneous \emph{prime} $D$-stables ideals
$\idp \subset \AnneauDePolynomes$ of rank $\geq n_1$.
Under these conditions there is a constant $K_1>0$ such that all $P \in
\AnneauDePolynomes \setminus \{0\}$ satisfying for all $C\geq C_1$
\begin{equation}\label{condition_n1D}
    i_0(\Z_C(P))\geq n_1
\end{equation}
satisfies also
\begin{equation} 
\ordz(P(\b{z},\ull{f})) \leq K_1(\deg_{\ul{X'}} P + 1)(\deg_{\ul{X}} P + 1)^n.
\end{equation}
\end{theorem}
\begin{proof} We are going to apply Theorem~\ref{LMGP}. The condition~(\ref{theoLMGP_condition_de_correctitude}) is satisfied because
the application $\phi=D$ is correct with respect to all the ideals (see Corollary~\ref{exemple_OpDiff_Correcte}).

The functions~$f_1(\b{z}),\dots,f_n(\b{z})$ are algebraically independent over $\kk(\b{z})$ because otherwise the ideal
$\idp_{\ullt{f}} \ne (0)$ been $D$-stable 
can not satisfy~(\ref{ordIleqKdegI2}).

Consider a bi-homogeneous equidimensional and $D$-stable ideal $J\subset\AnneauDePolynomes$, $J\ne(0)$. Let $\idq$ be a primary factor of $J$.
As $J$ is an equidimensional ideal, we can present $J$ as $J=\idq\cap I$, where $I\not\subset\sqrt{\idq}$. Then we fix an element $x\in I\setminus\sqrt{\idq}$.

Suppose $a\in\idq$, then $xa\in J$ and as $J$ is $D$-stable we have
\begin{equation}\label{Dxa_in_J}
    D(xa)\in J\subset \idq.
\end{equation}
Since $D$ is a derivation,
\begin{equation}\label{LMGPD_Leibnitz}
    D(xa)=D(x)a+xD(a).
\end{equation}
We have $D(x)a\in\idq$ (because $a\in\idq$), and also (given~(\ref{Dxa_in_J}) and~(\ref{LMGPD_Leibnitz})) $xD(a)\in\idq$.
By our choice of $x$ this implies $D(a)\in\idq$, thus $\idq$ is $D$-stable.

If $\car\kk=0$, we claim that $\sqrt{\idq}$ is $D$-stable. Indeed, for all $a\in\sqrt{\idq}$ we have $a^n\in\idq$, and so
\begin{equation*}
    D(a^n)\in\idq.
\end{equation*}
But $D(a^n)=nD(a^{n-1})$, thus in the case $\car\kk=0$ we deduce $D(a^{n-1})\idq$. By repeating this procedure we
arrive to the conclusion $D(a)\in\idq\subset\sqrt\idq$, thus $\sqrt{\idq}$ is $D$-stable.

Now the hypothesis~(\ref{ordIleqKdegI2}) assures immediately (\ref{RelMinN2}). Indeed,
\begin{eqnarray*}
    \dd_{(0,n-\rg\idp+1)}\sqrt{\idq}&\leq&\dd_{(0,n-\rg\idp+1)}\idq\leq\dd_{(0,n-\rg\idp+1)}\deg J,\\
    \dd_{(1,n-\rg\idp)}\sqrt{\idq}&\leq&\dd_{(1,n-\rg\idp)}\idq\leq\dd_{(1,n-\rg\idp)}\deg J
\end{eqnarray*}
and for at least one of the primary factors of $J$ we have
\begin{equation*}
    \ord_{\ullt{f}}\sqrt{\idq}=\ord_{\ullt{f}}\idq=\ord_{\ullt{f}}J.
\end{equation*}
So we can apply Theorem~\ref{LMGP}
(with the constant $C_0=0$, for example) and this theorem gives us the desired conclusion.
\end{proof}

\subsection{Zero order estimates for functions satisfying functional equations of generalized Mahler's type}


Let $\T$ be a rational  transformation from $\mpp^1 \times \mpp^{n}$ to itself
defined by
\begin{multline} \label{defT2}
(X_0':X_1',X_0:...:X_n)\rightarrow\Big(A_0'(X_0',X_1'):A_1'(X_0',X_1'),\\
A_0(X_0',X_1',X_0,...,X_n):...:A_n(X_0',X_1',X_0,...,X_n)\Big),
\end{multline}
where $A_i'\in\kk[X_0',X_1']$, $i=0,1$, are homogeneous polynomials of degree $r$ in $\ul{X}'$ and $A_j\in\AnneauDePolynomes$, $j=0,...,n$,
are bi-homogeneous polynomials of bi-degree $(s,t)$ in $\ul{X}'$ and $\ul{X}$.

\begin{remark}\label{rem_Mutual_Association} We associate to every rational transformation $\T$ defined by~(\ref{defT2}) and such that $A_0$, $A_0'$
are non-zero polynomials,
a system of functional equations~(\ref{relsTopfer}) by defining $p(\b{z})=\frac{A_1'(1,\b{z})}{A_0'(1,\b{z})}$:
\begin{equation} \label{relsTopfer2}
    A_0(\ullt{f}(\b{z}))f_i(p(\b{z}))=A_i(\ullt{f}(\b{z})), \quad i=1,...,n
\end{equation}
(where $\ullt{f}$ denotes $(1,\b{z},1,f_1(\b{z}),...,f_n(\b{z}))$). 

The other way around, starting from the system~(\ref{relsTopfer}) (where $p(\b{z})\in\kk(\b{z})$ and we do not make the hypothesis that $\ordz p \geq 2$)
the formulas~(\ref{defT2}) define a morphism $\T:\mpp^1 \times \mpp^{n}\rightarrow\mpp^1 \times \mpp^{n}$.
\end{remark}

\begin{definition} \label{def_Mutual_Association}
We say that the morphism defined by~(\ref{defT2}) and the system~(\ref{relsTopfer2}) are mutually associated.
\end{definition}

\begin{definition}\label{def_irrT}
The morphism $\T$ been defined by~(\ref{defT2}), we denote by $\irrT$ a union of zero locus of polynomial bi-homogeneous systems $A_i'(X_0',X_1',X_0,...,X_n)$, $i=0,1$, and $A_j(X_0',X_1',X_0,...,X_n)$, $j=0,...,n$. One has $\irrT \subset \mpp_{\kk}^1\times\mpp_{\kk}^n$ and this is a set of points where bi-projective application $\T$ is not defined (if $\irrT=\emptyset$ the transformation $\T$ is a regular bi-projective transformation).
\end{definition}

\begin{remark} \label{rem_TV_simplified}
To simplify the notation we write $\T(W)$ instead of $\T(W\setminus\irrT)$.
\end{remark}

\begin{definition}\label{definVarieteTstable}
We say that a sub-variety $W\subset\mpp^1\times\mpp^n$
is $\T$-stable, if
\begin{equation*}
 \ol{\T(W)}=W.
\end{equation*}
\end{definition}

\begin{remark} \label{varietestable_et_idealstable}
If a variety $W$ is $\T$-stable then the ideal $\I(W)$ is $\Talg$-stable, but the reciprocal statement is not true. The condition $\Talg(\I(W))\subset\I(W)$ geometrically means only that $W$ is a sub-scheme of $\T^{-1}(W)$. However, if we impose that the variety $W$ is irreducible and $\dim \T(W)=\dim W$, then
\begin{equation}
 \mbox{a variety $W$ is $\T$-stable} \Leftrightarrow \mbox{the ideal $\I(W)$ is $\Talg$-stable}.
\end{equation}
\end{remark}


We fix a point $\ullt{f}=(1:\b{f}'_1,1:\b{f}_1:...:\b{f}_n) \in
\mpp^1_{\KK}\times\mpp^n_{\KK}$ such that $\b{f}'_1$,$\b{f}_1$,...,$\b{f}_n$ are algebraically
independent over $\kk$.

\begin{lemma} \label{Testbon1} Let $\T:\mpp^1\times\mpp^n\rightarrow\mpp^1\times\mpp^n$ be a rational dominant transformation. There is a constant $\Cr$ depending on $\T$ and $\ullt{f}$ only, such that all irreducible variety
$V\subset\mpp^1\times\mpp^n$ satisfying
\begin{equation} \label{ordVpgC}
\ord_{\ullt{f}}V>\Cr,
\end{equation}
satisfies also
$$
\dim\T(V)=\dim V.
$$
\end{lemma}
\begin{proof}
We shall use the fact that the point
$\ullt{f}$ does not belong to any strict sub-variety of
$\mpp_{\kk}^1\times\mpp_{\kk}^n$ defined over $\kk$ (as $\b{f}'_1$,$\b{f}_1$,...,$\b{f}_n$ are algebraically
independent over $\kk$ by assumption).

Note that the morphism $\T$
is a regular one if restricted to some Zariski-open set of $\mpp_{\kk}^1\times\mpp_{\kk}^n$
defined over $\kk$, i.e. it is a regular morphism
apart from some Zariski-closed set $F_r$ (different from $\mpp_{\kk}^1\times\mpp_{\kk}^n$) defined
over $\kk$.

Then, apart from a closed set $F_e$ different from the whole space,
$\T$ is an \'etale morphism (proposition 4.5,
expos\'e~I \cite{SGA1}), hence it preserves the dimension of
varieties.

We define $F=F_r \cup F_e$. The variety $F$ is defined over
$\kk$ and is different from the whole space (because $F_r$ and $F_e$ are Zariski-closed strict subsets), we have $\ullt{f}\not\in F$, thus
$\dist(\ullt{f},F)=d_1>0$, or, alternatively,
$\ord(\ullt{f},F)=o_1 < +\infty$.

Let $V$ be an irreducible variety satisfying~(\ref{ordVpgC}) for
$\Cr=o_1$. By definition it signifies that at least one
point $\ull{\alpha} \in V$ satisfies
$\ord(\ullt{f},\ull{\alpha})>o_1$ and thus $V$ is not contained in $F$. So, in view of irreducibility of $V$, $U = V \setminus F$ is dense in $V$, hence $\dim V = \dim U$. As $\T$ is \'etale over $U$, we obtain $\dim \T U = \dim U$ and thus $\dim \T V \geq \dim \T U = \dim U = \dim V$. Inequality $\dim \T V \leq \dim V$ being true in any case, we conclude $\dim \T V = \dim V$.

So we arrive to the conclusion of the lemma with $\Cr=o_1$.
\end{proof}

Using~(\ref{LdTordZ}) or else
Theorem~\ref{dist_alpha}, we can assure that the cycle
$\b{Z}_{b}(P)$ contains at least one point $\alpha$ at a distance $b^{\frac{1}{n}}c_n^{-\frac{n-1}{n}}$ from $\ull{f}$, i.e. sufficiently close to $\ull{f}$. The next Lemma uses this remark in order to establish the alternative: either the polynomial $P$
satisfies a multiplicity lemma with some absolute constant
$C$, or we can suppose certain particular properties of action of $\T$ on any irreducible variety containing $\b{Z}_{C}(P)$.

\begin{lemma} \label{Testbon} {\it There is a constant $C_{\T}$
depending on $\T$ and $\ull{f}$ only, such that if a polynomial $P$
satisfies~(\ref{ordPplusqueb}) with $C \geq C_{\T}$, then $\T$
preserves the dimension of all irreducible varieties
containing $\b{Z}_{C}(P)$ (and so the transformation~$\Talg:\AnneauDePolynomes\rightarrow\AnneauDePolynomes$ is correct within the meaning of Definition~\ref{defin_phiestcorrecte}). }
\end{lemma}
\begin{proof} This is a corollary of Lemma~\ref{Testbon1}. The latter lemma gives us
a constant $\Cr$ (depending on $\T$).

Let $P$ be a polynomial satisfying~(\ref{ordPplusqueb}) with 
$C \geq c_n^{{n-1}}\Cr^{n}$. Using~(\ref{LdTordZ}) we find out that at least
one point $\ull{\alpha}$ of $\b{Z}_{C}(P)$
satisfies $\ord(\ull{f},\ull{\alpha})>\Cr$. So every irreducible
variety $W$ defined over $\kk(\b{z})$ and passing by $\ull{\alpha}$ satisfies
$\ord_{\ull{f}}W > \Cr$ and so, by Lemma~\ref{Testbon1},
$\T$ preserves the dimension of $W$.

Defining $C_{\T}=c_n^{{n-1}}\Cr^{n}$ we conclude.
\end{proof}

\begin{theorem}\label{LMGPF} Let $\kk$ be a field of an arbitrary characteristic and
$\T: \mpp^1_{\kk}\times\mpp^n_{\kk} \rightarrow \mpp^1_{\kk}\times\mpp^n_{\kk}$ a rational dominant transformation defined as in~(\ref{defT2}),
by the homogeneous polynomials $A_i'$, $i=0,1$ in $\ul{X}'$ of degree $r$, and polynomials $A_i$, $i=0,\dots,n$ bi-homogeneous in $\ul{X}'$ and $\ul{X}$,
of bi-degree $(s,t)$.

Let $f_1(\b{z})$,...,$f_n(\b{z}) \in \kk[[\b{z}]]$ be power series algebraically independent over $\kk(\b{z})$
and $n_1\in\{1,\dots,n\}$, $C_1\in\mrr^+$.
Wa also denote $\ullt{f}=(1,\b{z},1,f_1(\b{z}),...,f_n(\b{z}))$.

Suppose moreover that there exists $\lambda \in \mrr_{>0}$, such that for all $Q\in\AnneauDePolynomes$
\begin{equation} \label{section_AB_IElambda}
\Ordz Q(\T(\ullt{f})) \geq \lambda \Ordz Q(\ullt{f}),
\end{equation}
and that exists a
constant $K_0 \in \mrr^{+}$ (dependent on $\T$ and
$\ullt{f}$ only) such that for all positive integer
\begin{equation}\label{Nmajoration}
    N\leq\nu(n+1)!\rho_{n+1}^{n+1}
\end{equation}
(where $\nu$ equals $2^{n+2}\max\left(1,\frac{4r}{s}\right)^{n+1}$ if $s\ne 0$ and 1 otherwise)
all irreducible ${\T}^N$-stable variety $W\varsubsetneq\mpp^1\times\mpp^n$ (defined over the field $\kk$) of dimension $\dim W\leq n-n_1+1$
satisfies necessarily
\begin{equation} \label{RelMinN}
\ord_{\bt{f}}(W) < K_0\left(\dd_{(0,\dim W)}W+\dd_{(1,\dim W-1)}W\right).
\end{equation}
Then there exists a constant $K_1>0$ such that for all $P \in \AnneauDePolynomes\setminus\{0\}$ satisfying for all $C\geq C_1$
\begin{equation}\label{condition_n1F}
    i_0(\Z_C(P))\geq n_1,
\end{equation}
satisfies also
\begin{equation} \label{LdMpolynome}
\ordz(P(\ullt{f})) \leq K_1(\deg_{\ul{X'}}P + \deg_{\ul{X}}P + 1)(\deg_{\ul{X}} P + 1)^n.
\end{equation}
\end{theorem}
\begin{proof} We are going to apply Theorem~\ref{LMGP} to the transformation $\phi=\Talg$.
One readily verify properties~(\ref{degphiQleqdegQ}) (with $\mu=t$, $\nu_0=r$ and $\nu_1=s$) and~(\ref{condition_T2_facile})
(the latter is provided by hypothesis~(\ref{section_AB_IElambda})). In this case
$\nu$ (introduced in Definition~\ref{def_nu}) is equal to $2^{n+2}\max\left(1,\frac{4r}{s}\right)^{n+1}$ if $s\ne 0$ and $1$ otherwise.
Define $C_0=\Cr$. Then the condition~(\ref{theoLMGP_condition_de_correctitude}) is assured
by Lemma~\ref{Testbon1} and Corollary~\ref{exemple_ApBirr_Correcte}.

In order to verify~(\ref{RelMinN2}) for equidimensional $\Talg$-stable ideals $J$
satisfying $m(J)\leq\nu(n+1)\rho_{n+1}^{n+1}$ and~(\ref{theoLMGP_condition_ordp_geq_C0}),
we note that the condition of $\Talg$-stability implies
\begin{equation}\label{VJsubsetTm1VJ}
    \V(J)\subset\ol{\T^{-1}(\V(J))}
\end{equation}
whence by application of $\T$ to the both sides
\begin{equation}\label{TVJsubsetVJ}
    \T(\V(J)\setminus\irrT)\subset\V(J)
\end{equation}
(we recall that $\irrT$ is introduced in Definition~\ref{def_irrT}). As we have defined $C_0=\Cr$ (and in view of~(\ref{theoLMGP_condition_ordp_geq_C0})
and Lemma~\ref{Testbon1})
each component $W$ of variety $\V(J)$ satisfies $W\not\subset\irrT$
\begin{equation*}
    \dim W=\dim\T(W\setminus\irrT).
\end{equation*}
We recall that $J$ is an equidimensional ideal, hence for all its irreducible components one has$\dim W=\dim\V(J)$, thus
$\ol{\T(W\setminus\irrT)}$
is an irreducible component of $\V(J)$.

So $\T$ induces an application from
the set of irreducible components of $\V(J)$ to itself. But the cardinality of the set of irreducible components of $\V(J)$ is estimated from above by
$m(J)\leq\nu(n+1)!\rho_{n+1}^{n+1}$ by hypothesis on $J$.
Thus there is a subset of this set where $\T$ acts as a permutation, and this subset
contains at most $\nu(n+1)!\rho_{n+1}^{n+1}$ elements. Orbits of $\T$ in this subset can not have more than
$\nu(n+1)!\rho_{n+1}^{n+1}$ elements, so we deduce that some iteration $\T^N$,
$N\leq\nu(n+1)!\rho_{n+1}^{n+1}$, fixes an irreducible component $W$
of $\V(J)$.

Now the hypothesis~(\ref{RelMinN}) is assured by hypothesis~(\ref{RelMinN2}), and also
all the conditions of Theorem~\ref{LMGP} are verified. Applying this theorem
we obtain conclusion~(\ref{LdMpolynome}) (this conclusion coincides with conclusion~(\ref{LdMpolynome2})
of Theorem~\ref{LMGP} if $K_1=K\mu^{n-1}\max(\mu+\nu_0,\nu_1)$).
\end{proof}

\subsection{Linear case} \label{sectionNishioka}

In this subsection we let $\T:\mpp^1\times\mpp^n\rightarrow\mpp^1\times\mpp^n$ denote a transformation defined by~(\ref{defT2}) with polynomials $A_i$ linear in $\ul{X}$ (and so $t=1$):
\begin{equation} \label{defAlinear}
    A_i=\sum_{j=0}^n a_{ij}(\ul{X}')X_j, \qquad i=0,...,n,
\end{equation}
where $a_{ij}(\ul{X}')$ denote homogeneous polynomials from $\kk[X_0',X_1']$ of degree $s$. Suppose that system~(\ref{defT2}) in this situation admits a solution in algebraically independent power series, we denote this solution by $(1,f_1(\b{z}),...,f_n(\b{z}))$. One readily verifies that in this situation the matrix $(a_{ij})_{i,j=1,\dots,n}$ is non-degenerated, therefore the transformation $\T$ is dominant and $\irrT=\emptyset$ (the proof of a more general statement can be found in the beginning of chapter~3 of~\cite{EZ}, see Lemma~3.3 \emph{Op. cit.}).

In this situation we can deduce from Theorem~\ref{LMGPF} an unconditional optimal multiplicity lemma, Theorem~\ref{theoNishioka}. Its proof is the matter of this section. Our unconditional theorem can be considered as a counterpart of the similar result for systems of linear differential equations due to Yu.Nesterenko~\cite{N1974}.

\smallskip

From now on and up to the end of this section we denote by $W$ a variety defined over $\kk$ and $\T$-stable.

\smallskip

Let $\gamma(\b{z})\in\KK$. We define
\begin{equation}\label{defWz}
    W_{\gamma(\b{z})}\eqdef W\cap\V(X_0'\gamma(\b{z})-X_1')
\end{equation}
and
$\T_{\gamma(\b{z})}$ the transformation from $\mpp^{n}_{\KK}$ to
$\mpp^{n}_{\KK}$ defined by
\begin{equation*}
\begin{aligned}
(X_0:...:X_n)\rightarrow\Big(A_0(1,\gamma(\b{z}),X_0,...,X_n):...:A_n(1,\gamma(\b{z}),X_0,...,X_n)\Big).
\end{aligned}
\end{equation*}


\smallskip

\begin{lemma} \label{sectionNishioka_teclem1} Let $\gamma(\b{z})\in\KK$ et $\bt{\phi}\in W_{\gamma(\b{z})}$.
We suppose
\begin{equation}\label{gamma_est_conj_p}
    p(\gamma(\b{z}))=p(\b{z}).
\end{equation}
Then there is a $\bt{\phi'}\in W_{\b{z}}$ such that $\T(\bt{\phi'})=\T(\bt{\phi})$.
\end{lemma}
\begin{proof}
Applications $\T_{\gamma(\b{z})}:\pi_1(W_{\gamma(\b{z})})\rightarrow\pi_1(W_{p(\b{z})})$ and $\T_{\b{z}}:\pi_1(W_{\b{z}})\rightarrow\pi_1(W_{p(\b{z})})$ both are linear and non-degenerated (because $\b{z}$ and $\gamma(\b{z})$ both are transcendental over $\kk$, while $\T$ is defined over $\kk$ and dominant. So these two applications induce isomorphisms between $\pi_1(W_{\gamma(\b{z})})$, $\pi_1(W_{p(\b{z})})$ and $\pi_1(W_{\b{z}})$. Now we readily verify that if we define $\b{\phi'}:=\T_{\b{z}}^{-1}\circ\T_{\gamma(\b{z})}(\pi_1(\bt{\phi}))$, the point $\bt{\phi'}:=(1:\b{z},\b{\phi'})$ belongs to $W_{\b{z}}$ and we have $\T(\bt{\phi'})=\T(\bt{\phi})$.
\end{proof}


\begin{lemma} \label{sectionNishioka_teclem2}
Let $\bt{\beta}=(1:\b{z},\b{\beta}(\b{z}))\in W$. Then $W\cap\T^{-1}(1:p(\b{z}),\b{\beta}(p(\b{z})))\ne\emptyset$.
\end{lemma}
\begin{proof}
As $W$ is defined over $\kk$ and $\bt{\beta}(\b{z})\in W$, we have $\bt{\beta}(p(\b{z}))\in W$. As $W$ is $\T$-stable, $W\cap\T^{-1}(\bt{\beta}(p(\b{z})))\ne\emptyset$.
\end{proof}


\begin{corollary} \label{sectionNishioka_teclem3}
Let $\bt{\beta}=(1:\b{z},\b{\beta}(\b{z}))\in W_{\b{z}}$. Then $W_{\b{z}}\cap\T^{-1}(1:p(\b{z}),\b{\beta}(p(\b{z})))\ne\emptyset$.
\end{corollary}
\begin{proof}
The point $\bt{\beta}$ being given, we find with Lemma~\ref{sectionNishioka_teclem2} a point $\bt{\alpha}$ in $W\cap\T^{-1}(\bt{\beta}(p(\b{z})))$. As $\T(\bt{\alpha})=\bt{\beta}(p(\b{z}))\in W_{p(\b{z})}$, we necessarily have $\bt{\alpha}=(1:\gamma(\b{z}),\b{\alpha}(\b{z}))$, where $\gamma(\b{z})$ satisfies $p(\gamma(\b{z}))=p(\b{z})$ (and thus by definition $\bt{\alpha}\in W_{\gamma(\b{z})}$). Hence we can apply Lemma~\ref{sectionNishioka_teclem1}, which provides a point $\bt{\alpha}'\in W_{\b{z}}$ satisfying $\T(\bt{\alpha}')=\T(\bt{\alpha})=\bt{\beta}(p(\b{z}))$ and this completes the proof.
\end{proof}

\begin{lemma} \label{sectionNishioka_teclem4}
\begin{equation} \label{sectionNishioka_ieq1}
    \sup_{\alpha\in\pi_1(W_{\b{z}})}\Ordz(\Tbz(\b{f}),\Tbz(\b{\alpha}))\geq\lambda\sup_{\beta\in\pi_1(W_{\b{z}})}\Ordz(\b{f},\b{\beta}).
\end{equation}
\end{lemma}
\begin{proof}
We note that
\begin{equation} \label{sectionNishioka_ieq2}
    \lambda\sup_{\beta\in\pi_1(W_{\b{z}})}\Ordz(\b{f},\b{\beta})=\sup_{\beta\in\pi_1(W_{\b{z}})}\Ordz(\b{f}(p(\b{z})),\b{\beta}(p(\b{z}))).
\end{equation}

For $\ull{f}$ we have the property $\ull{f}(p(\b{z}))=\Tbz(\ull{f})$. Then for all $\b{\beta}\in\pi_1(W_{\b{z}})$ there is a point $\bt{\beta}\in W_{\b{z}}$ satisfying $\pi_1(\bt{\beta})=\b{\beta}$ and then Corollary~\ref{sectionNishioka_teclem3} gives us a point $\bt{\alpha}'\in W_{\b{z}}$ satisfying $\T(\bt{\alpha}')=\bt{\beta}(p(\b{z}))$, hence
\begin{equation} \label{sectionNishioka_alphabeta}
    \Tbz(\pi_1(\bt{\alpha}'))=\b{\beta}(p(\b{z})).
\end{equation}
So for every point $\b{\beta}\in\pi_1(W_{\b{z}})$ we have determined a point $\b{\alpha}=\pi_1(\bt{\alpha}')\in \pi_1\left(W_{\b{z}}\right)$ satisfying~(\ref{sectionNishioka_alphabeta}).

This construction, coupled with~(\ref{sectionNishioka_ieq2}), readily gives us~(\ref{sectionNishioka_ieq1}).
\end{proof}

We are ready to prove the principal result of this section.


\begin{proof}[Proof of Theorem~\ref{theoNishioka}]
 Note that $\T$ is defined by~(\ref{defT2}) with all the polynomials $A_i$ linear in $\ul{X}$. Thus the transformation $\T^N$ is also of this type. Further, as $N$ has an upper bound that depends on $\T$ only, the corresponding degrees $s$, $\lambda$ etc. of the transformation $\T^N$ also have upper bounds that depend on $\T$ only. So we bring the analysis of ideals stable under the transformation $\T^N$ to the analysis of $\T$-stables ideals.

Note that the application $\Tbz$ is a linear non-degenerated application defined over $\kk(\b{z})$. We introduce one more notation: $C_3\eqdef\ordz\det\Tbz^{-1}$.

Let $W\subset\mpp^1\times\mpp^n$ be a $\T$-stable variety defined over $\kk$. We claim that
\begin{equation} \label{theoremeNishioka_ieq1}
    \sup_{\b{\alpha}\in\pi_1(W_{\b{z}})}\Ordz(\b{f},\b{\alpha})\leq\frac{C_3}{\lambda-1}.
\end{equation}
Indeed, by Lemma~\ref{sectionNishioka_teclem4} we have~(\ref{sectionNishioka_ieq1}). Further, as $\Tbz$ is a linear non-degenerated application, one has $\Ordz(\Tbz(\b{f}),\Tbz(\b{\alpha}))\leq C_3+\Ordz(\b{f},\b{\alpha})$ and thus
\begin{equation} \label{sectionNishioka_ieq1_1}
    \sup_{\b{\alpha}\in\pi_1(W_{\b{z}})}\Ordz(\Tbz(\b{f}),\Tbz(\b{\alpha}))\leq C_3+\sup_{\b{\alpha}\in\pi_1(W_{\b{z}})}\Ordz(\b{f},\b{\alpha}).
\end{equation}
Comparing~(\ref{sectionNishioka_ieq1}) and~(\ref{sectionNishioka_ieq1_1}) we find~(\ref{theoremeNishioka_ieq1}).

We claim that inequality~(\ref{theoremeNishioka_ieq1}) implies
\begin{equation} \label{theoremeNishioka_ieq1_2}
    \ord_{\ullt{f}}W=\sup_{\bt{\alpha}\in W}\Ordz(\bt{f},\bt{\alpha})\leq \max\left(2,\frac{C_3}{\lambda-1}\right).
\end{equation}
Indeed, if for an $\bt{\alpha}\in W$ one has
\begin{equation}\label{theoremeNishioka_Ord_tildealpha_f_est_grand}
    \Ordz(\bt{f},\bt{\alpha})>\max\left(2,\frac{C_3}{\lambda-1}\right),
\end{equation}
then there is a system of projective coordinates $\ullt{\alpha}=(1,\gamma(\b{z}),\ull{\alpha}(\b{z}))$ with $\gamma(\b{z})\in\ol{\kk((\b{z}))}$ and $\b{\alpha}(\b{z})\in\pi_1(W)$ such that $\ordz(\b{z}-\gamma(\b{z}))>\max(2,\frac{C_3}{\lambda-1})$ and
\begin{equation}\label{theoremeNishioka_Ord_alpha_f_est_grand}
    \Ordz(\b{\alpha}(\b{z}),\b{f})>\max(2,\frac{C_3}{\lambda-1}).
\end{equation}
As $\ordz(\b{z}-\gamma(\b{z}))>\max(2,\frac{C_3}{\lambda-1})$, there exists a series $v(\b{z})\in\kk((\b{z}))$ satisfying $\gamma(v(\b{z}))=\b{z}$ and
\begin{equation}\label{ord_v_grand}
    \ordz(\b{z}-v(\b{z}))>\max\left(2,\frac{C_3}{\lambda-1}\right).
\end{equation}
Substituting $v(\b{z})$ to $\ullt{\alpha}(\b{z})$ in place of $\b{z}$ we find $\ullt{\alpha}(v(\b{z}))=(1,\b{z},\b{\alpha}(v(\b{z})))$. Variety $W$ being defined over $\kk$, one has $\ullt{\alpha}(v(\b{z}))\in W$.  Moreover, as $v(\b{z})$ coincides with the series $\b{z}$ at least to the order $\max(2,\frac{C_3}{\lambda-1})$, inequality~(\ref{theoremeNishioka_Ord_alpha_f_est_grand}) implies
\begin{equation*} 
    \Ordz(\b{\alpha}(v(\b{z})),\b{f})>\max\left(2,\frac{C_3}{\lambda-1}\right).
\end{equation*}
As $\ullt{\alpha}\in W_{\b{z}}$, we deduce the lower bound
\begin{equation*}
    \sup_{\b{\alpha}\in\pi_1(W_{\b{z}})}\Ordz(\b{f},\b{\alpha})>\max\left(2,\frac{C_3}{\lambda-1}\right),
\end{equation*}
 that contradicts~(\ref{theoremeNishioka_ieq1}). This contradiction shows that assumption~(\ref{theoremeNishioka_Ord_tildealpha_f_est_grand}) is untenable and therefore assures~(\ref{theoremeNishioka_ieq1_2}).

Now we remark that the transformation $\T^N$, $N\in\mnn$, is also defined by formulas of type~(\ref{defT2}), so any $\T^N$-stable variety
satisfies~(\ref{theoremeNishioka_ieq1_2}) (with $C_3$ replaced by $C_3^N$).
The latter inequality shows that hypothesis~(\ref{RelMinN}) of Theorem~\ref{LMGPF} is satisfied (even in a stronger form).
It is sufficient now to apply Theorem~\ref{LMGPF} to obtain the desired result.
\end{proof}


\begin{center}%
          {\bfseries Acknowledgement\vspace{-.5em}}%
\end{center}%
\thanks{The author would like to express his profound gratitude to Patrice \textsc{Philippon}. His interventions at many stages of  this research was of decisive importance.}

{\small



\def\cprime{$'$} \def\cprime{$'$} \def\cprime{$'$} \def\cprime{$'$}
  \def\cprime{$'$} \def\cprime{$'$} \def\cprime{$'$} \def\cprime{$'$}
  \def\cprime{$'$} \def\cprime{$'$} \def\cprime{$'$} \def\cprime{$'$}
  \def\cprime{$'$} \def\cprime{$'$} \def\cprime{$'$} \def\cprime{$'$}
  \def\polhk#1{\setbox0=\hbox{#1}{\ooalign{\hidewidth
  \lower1.5ex\hbox{`}\hidewidth\crcr\unhbox0}}} \def\cprime{$'$}
  \def\cprime{$'$}

\bigskip

\noindent{\footnotesize EZ\,: }\begin{minipage}[t]{0.9\textwidth}
\footnotesize{\sc Institut de math\'ematiques de Jussieu, Universit\'e Paris 7, Paris, France}\\
{\it E-mail address}\,:~~ \verb|zorin@math.jussieu.fr|\quad\emph{or}\quad\verb|EvgeniyZorin@yandex.ru|
\end{minipage}

}

\end{document}